\documentclass[12pt]{article}
\usepackage{amsfonts}
\usepackage{amsthm}
\usepackage{amsmath}
\usepackage{amssymb}
\usepackage{latexsym}
\usepackage{bbm}
\usepackage[pdftex]{graphicx,color}
\usepackage[table,xcdraw]{xcolor}
\usepackage{enumerate}
\usepackage{verbatim}
\usepackage{amsfonts}
\usepackage{epsfig}
\usepackage{tikz}

\usepackage{colortbl}
\definecolor{text1}{cmyk}{1,.65,0,0} 
\definecolor{text2}{rgb}{1,0,0} 
\definecolor{text3}{cmyk}{0,0,0,1} 
\definecolor{text4}{cmyk}{0,0,0,0.5} 

\usepackage{bbm} 
\usepackage{amssymb,amsmath,tikz,graphicx,float,xspace,chemarrow}
\usepackage{graphicx}
\usepackage[utf8]{inputenc}
\makeatletter

\newcommand{\be}{\begin{equation}}
\newcommand{\ee}{\end{equation}}
\newcommand{\bea}{\begin{eqnarray}}
\newcommand{\eea}{\end{eqnarray}}
  
\newcommand{\distas}[1]{\mathbin{\overset{#1}{\kern\z@\sim}}}%
\newsavebox{\mybox}\newsavebox{\mysim}
\newcommand{\distras}[1]{%
  \savebox{\mybox}{\hbox{\kern3pt$\scriptstyle#1$\kern3pt}}%
  \savebox{\mysim}{\hbox{$\sim$}}%
  \mathbin{\overset{#1}{\kern\z@\resizebox{\wd\mybox}{\ht\mysim}{$\sim$}}}%
}
\newcommand{\norm}[1]{
\|{#1}\|^2}
\definecolor{mygray}{gray}{0.38}
\definecolor{davysgrey}{rgb}{0.33, 0.33, 0.33}

\newtheorem{theorem}{Theorem}[section]
\newtheorem{lemma}{Lemma}[section]
\newtheorem{remark}{Remark}[section]
\newtheorem{proposition}{Proposition}[section]
\newtheorem{corollary}{Corollary}[section]
\newtheorem{example}{Example}[section]
\newtheorem{definition}{Definition}[section]

\makeatletter

\@addtoreset{equation}{section}
\makeatother

\begin{document}

\begin{center}
   {\bf On predictive density estimation with additional information \footnote{\today}}  \end{center}

\begin{center}
{\sc \'Eric MARCHAND$^{a}$, Abdolnasser SADEGHKHANI$^{b}$} \\

{\it a  Universit\'e de
    Sherbrooke, D\'epartement de math\'ematiques, Sherbrooke (Qu\'{e}bec),
    CANADA \quad (eric.marchand@usherbrooke.ca) } \\
    
{\it b  Queen's University, Department of Mathematics and Statistics, Kingston (Ontario),     CANADA \quad (a.sadeghkhani@queensu.ca) } \\    
\end{center}
\vspace*{0.2cm}
\begin{center}
{\sc Summary} \\
\end{center}
\vspace*{0.1cm}
\small
Based on independently distributed $X_1 \sim N_p(\theta_1, \sigma^2_1 I_p)$ and  $X_2 \sim N_p(\theta_2, \sigma^2_2 I_p)$, we consider the efficiency of various predictive density estimators for $Y_1 \sim N_p(\theta_1, \sigma^2_Y I_p)$, with the additional information $\theta_1 - \theta_2 \in A$ and known $\sigma^2_1, \sigma^2_2, \sigma^2_Y$.  We provide improvements on benchmark predictive densities such as {\em plug-in}, the maximum likelihood, and the minimum risk equivariant predictive densities.   Dominance results are obtained for $\alpha-$divergence losses and include Bayesian improvements for reverse Kullback-Leibler loss, and Kullback-Leibler (KL) loss in the univariate case ($p=1$).    An ensemble of techniques are exploited, including variance expansion (for KL loss),  point estimation duality, and concave inequalities.  Representations for Bayesian predictive densities, and in particular for  $\hat{q}_{\pi_{U,A}}$ associated with a uniform prior for $\theta=(\theta_1, \theta_2)$ truncated to $\{\theta \in \mathbb{R}^{2p}: \theta_1 - \theta_2 \in A  \}$, are established and are used for the Bayesian dominance findings.  Finally and interestingly, these Bayesian predictive densities also relate to skew-normal distributions, as well as new forms of such distributions.

\normalsize

\vspace*{0.5cm}


\normalsize
\section{Introduction}
\label{intro}
\subsection{Problem and Model}
Consider independently distributed
\begin{align}
\label{model}
X=\displaystyle \binom{X_1}{X_2} \sim\mathrm{N}_{2p}\left(\displaystyle \theta=\binom{\theta_1}{\theta_2},\,\Sigma=\bigl( \begin{smallmatrix} 
  \sigma_1^2 I_{p} & 0\\
  0 & \sigma_2^2 I_{p} 
\end{smallmatrix} \bigr)\right)\,,   Y_1 \sim\mathrm{N}_{p}(\theta_1,\,\sigma^{2}_Y I_{p})\,,
\end{align}
where $X_1,X_2, \theta_1, \theta_2$ are $p-$dimensional, and
with the additional information (or constraint) $\theta_1-\theta_2 \in A \subset \mathbb{R}^p$, $A $, $\sigma_1^{2}, \sigma_2^{2}$, $\sigma^{2}_Y$ all known, the variances not necessarily equal.  We investigate how to gain from the additional information in providing a predictive density $\hat{q}(\cdot;X)$ as an estimate of the density $q_{\theta_1}(\cdot)$ of $Y_1$.  Such a density is of interest as a surrogate for $q_{\theta_1}$, as well as for generating either future or missing values of $Y_1$.  The additional information $\theta_1 - \theta_2 \in A$ renders $X_2$ useful in estimating the density of $Y_1$ despite the independence and the otherwise unrelated parameters.

The reduced $X$ data of the above model is pertinent to summaries $X_1$ and $X_2$ that arise through a sufficiency reduction, a large sample approximation, or limit theorems.   Specific forms of $A$ include:
\begin{enumerate}
\item[ {\bf (i)}]
order constraints $\theta_{1,i}-\theta_{2,i} \geq 0$ for $i=1,\ldots, p\,$; the $\theta_{1,i}$ and $\theta_{2,i}$'s representing the components of $\theta_1$ and $\theta_2$;

\item[ {\bf (ii)}]
rectangular constraints $|\theta_{1,i} - \theta_{2,i}| \leq m_i$ for $i=1,\ldots, p\,$;

\item[ {\bf (iii)}] spherical constraints $\|\theta_1 - \theta_2\| \leq m \,$;
\item[ {\bf (iv)}] order and bounded constraints $m_1 \geq \theta_{1,i} \geq \theta_{2,i} \geq m_2$ for $i=1,\ldots, p\,$.  
\end{enumerate}

There is a very large literature on statistical inference in the presence of such constraints, mostly for {\bf (i)} (e.g.,  Hwang and Peddada, 1994;  Dunson and Neelon, 2003; Park, Kalbfleisch and Taylor, 2014) among many others).  Other sources on estimation in restricted parameter spaces can be found in the review paper of Marchand and Strawderman (2004), as well as the monograph by van Eeden (2006).  There exist various findings for estimation problems with additional information, dating back to Blumenthal and Cohen (1968) and  Cohen and Sackrowitz (1970), with further contributions by van Eeden and Zidek (2001, 2003), Marchand et al. (2012), Marchand and Strawderman (2004).

\begin{remark}
\label{transformation}
Our set-up applies to various other situations that can be transformed or reduced to model (\ref{model}) with $\theta_1 - \theta_2 \in A$.  Here are some examples.  

\begin{enumerate}
\item[ {\bf (I)}]

Consider model (\ref{model}) with the linear constrained $c_1\theta_1 - c_2 \theta_2 + d \in A$, $c_1, c_2$ being constants not equal 
to $0$, and $d \in \mathbb{R}^p$.  Transforming  $X_1'=c_1 X_1, X_2'=c_2X_2 - d$, and $Y_1'=c_1 Y_1$ leads to model (\ref{model}) based on the triplet $(X_1', X_2', Y_1')$, expectation parameters $\theta_1'=c_1\theta_1, \theta_2'=c_2\theta-d$, covariance matrices $c_i^2 \sigma^2_i I_p, i=1,2$ and $c_1^2 \sigma_Y^2 I_p$, and with the additional information $\theta_1' -\theta_2' \in A $.  With the class of losses being intrinsic (see Remark \ref{intrinsic}), and the study of predictive density estimation for $Y_1'$ equivalent to that for $Y_1$, our basic model and the findings below in this paper will indeed apply for linear constrained $c_1\theta_1 - c_2 \theta_2 + d \in A$.

\item[ {\bf (II)}]
Consider a bivariate normal model for $X$ with means $\theta_1, \theta_2$, variances $\sigma_1^2$, $\sigma_2^2$, correlation coefficient $\rho >0$, and the additional information $\theta_1 - \theta_2 \in A$.   The transformation  $X_1'=X_1$, $X_2'= \frac{1}{\sqrt{1+\rho^2}
} (X_2 - \frac{\rho \sigma_2}{\sigma_1} X_1)$ leads to independent coordinates with means $\theta_1'=\theta_1, \theta_2'=\frac{1}{\sqrt{1+\rho^2}
} (\theta_2 - \frac{\rho \sigma_2}{\sigma_1} \theta_1)$, and variances $\sigma_1^2$, $\sigma_2^2$.  We thus obtain model (\ref{model}) for $(X_1', X_2')$ with the additional information $\theta_1 - \theta_2 \in A$ transformed to  $c_1\theta_1' - c_2 \theta_2' + d \in A$, as in part {\bf (I)} above, with $c_1= 1 + \frac{\rho \sigma_2}{\sigma_1}$, $c_2= \sqrt{1+\rho^2}$, and $d=0$.

\end{enumerate}
\end{remark} 
 
\subsection{Predictive density estimation} 
\noindent Several loss functions are at our disposal to measure the efficiency of estimate $\hat{q}(\cdot;x)$, and these include the class of $\alpha-$divergence loss functions (e.g., Csisz\`{a}r, 1967) given by
\begin{equation}
\label{alphadivergence}
L_{\alpha}(\theta, \hat{q}) \,=\, \int_{\mathbb{R}^p} \, h_{\alpha}\left(\frac{\hat{q}(y;x)}{q_{\theta_1}(y)} \right) \; q_{\theta_1}(y)\, dy\,,
\end{equation} 
with 
 \[
h_{\alpha}(z)=
    \left\{
          \begin{array}{ll}
               {4\over 1-\alpha^2}(1-z^{(1+\alpha)/2})~~~\mbox{for}~ |\alpha|< 1\\
               z \log(z)~~~~~~~~~~~~~~~~\mbox{for}~ \alpha= 1\\
                -\log(z)~~~~~~~~~~~~~~~\mbox{for}~ \alpha= -1.
          \end{array}
    \right.
\] 
Notable examples in this class include Kullback-Leibler ($h_{-1}$), reverse Kullback-Leibler ($h_1$), and Hellinger ($h_0/4$).  The cases $|\alpha|<1$ stand apart, and merit study, in the sense that these losses are typically bounded, whereas both Kullback-Leibler and reverse Kullback-Leibler losses are typically unbounded (see Remark \ref{lim}).  
For an above given loss,   we measure the performance of a predictive density $\hat{q}(\cdot;X)$ by the frequentist risk 

\begin{equation}
\label{frequentistrisk}
R_{\alpha}(\theta, \hat{q}) \,=\,  \int_{\mathbb{R}^{2p}}  L_{\alpha}\left(\theta, \hat{q}(\cdot;x)\right) \,  \,p_{\theta}(x) \, dx  \,,  
\end{equation} 
$p_{\theta}$ representing the density of $X$.

Such a predictive density estimation framework was outlined for Kullback-Leibler loss in the pioneering work of Aitchison and Dunsmore (1975), as well as Aitchison (1975), and has found its way in many different fields of statistical science such as decision theory, information theory, econometrics, machine learning, image processing, and mathematical finance. 
There has been much recent Bayesian and decision theory analysis of predictive density estimators, in particular for multivariate normal or spherically symmetric settings, as witnessed by the work of Komaki (2001), George, Liang and Xu (2006), Brown, George and Xu (2008), Kato (2009), Fourdrinier et al. (2011), Ghosh, Mergel and Datta (2008), Maruyama and Strawderman (2012), Kubokawa, Marchand and Strawderman (2015, 2017), among others.

\begin{remark}
\label{intrinsic}
We point out that losses in (\ref{alphadivergence}) are intrinsic
in the sense that predictive density estimates of the density of 
$Y'=g(Y)$, with invertible $g: \mathbb{R}^p \to \mathbb{R}^p$ and inverse jacobian $J$, lead to an equivalent loss with the natural choice $\hat{q}(g^{-1}(y');x) \, |J|$ as
$$  
\int_{\mathbb{R}^p} \, h_{\alpha}\left(\frac{\hat{q}(g^{-1}(y');x) \, |J|}{q_{\theta_1}(g^{-1}(y')) \, |J|} \right) \; q_{\theta_1}(g^{-1}(y'))\, |J| \, dy'\,
=   \int_{\mathbb{R}^p} \, h_{\alpha}\left(\frac{\hat{q}(y;x)}{q_{\theta_1}(y)} \right) \; q_{\theta_1}(y)\, dy\,, $$
which is indeed $L_{\alpha}(\theta, \hat{q})$ independently of $g$.

\end{remark}

\subsection{Description of main findings}

In our predictive density estimation framework, we study various predictive densities such as: {\bf (i)} {\em plug-in} densities $N_p(\hat{\theta}_1(X), \sigma^2_Y I_p)$ including the predictive maximum likelihood estimator (mle); {\bf (ii)} minimum risk equivariant (MRE) predictive densities $\hat{q}_{mre}$; {\bf (iii)} variance expansions $N_p(\hat{\theta}_1(X), c \sigma^2_Y I_p)$, with $c>1$, of {\em plug-in} predictive densities; and {\bf (iv)}  Bayesian predictive densities with an emphasis on the uniform prior for $\theta$ truncated to the information set $A$. Our findings concern, except for Section 2, frequentist risk performance as in (\ref{frequentistrisk}), and related dominated dominance results covering the class of $\alpha-$divergence losses $L_{\alpha}$, as well as various types of information sets $A$.  

Subsection 3.1 provides Kullback-Leibler improvements on {\em plug-in} densities by variance expansion.  We make use of a technique due to Fourdrinier et al. (2011), which is universal with respect to $p$ and $A$ and requiring a determination, or lower-bound, of the infimum mean squared error of the plug-in estimator.  Such a determination is facilitated by a mean squared error decomposition (Lemma \ref{decomposition}) expressing the risk in terms of the risk of a one-population restricted parameter space estimation problem.  Such a decomposition appears in Marchand and Strawderman (2004).  

The dominance results of Subsection 3.2 apply to $L_{\alpha}$ losses and exploit point estimation duality.  The targeted predictive densities to be improved upon include {\em plug-in} densities, $\hat{q}_{mre}$, and more generally predictive densities of the form $\hat{q}_{\hat{\theta}_1, c} \sim N_p(\hat{\theta}_1(X), c \sigma^2_Y I_p)$.  The focus here is on improving on {\em plug-in} estimates $\hat{\theta}_1(X)$ by exploiting a correspondence with the problem of estimating $\theta_1$ under a dual loss.  Both Kullback-Leibler and reverse Kullback-Leibler losses lead to dual mean squared error performance.  In turn, as in Marchand and Strawderman (2004), the above risk decomposition relates this performance to a restricted parameter space problem.  Results for such problems are thus borrowable to infer dominance results for the original predictive density estimation problem.   For other $\alpha-$divergence losses, the strategy is similar, with the added difficulty that the dual loss relates to a reflected normal loss.  But, this is handled through a concave inequality technique (e.g., Kubokawa, Marchand and Strawderman, 2015) relating risk comparisons to mean squared error comparisons.  Several examples complement the presentation of Section 3.

Sections 2, 4, and 5 relate to Bayesian predictive densities, and especially to the Bayes procedure $\hat{q}_{\pi_{U,A}}$ with respect to the uniform prior $\mathbb{I}_A(\theta_1-\theta_2)$ restricted to $A$.   Section 2 presents various representations for $\hat{q}_{\pi_{U,A}}$, with examples connecting not only to known skewed-normal distributions, but also to seemingly new families of skewed-normal type distributions.  Section 4 contains Bayesian dominance results for both reverse Kullback-Leibler and Kullback-Leibler losses.  The case of reverse Kullback-Leibler loss, which is addressed in Subsection 4.1, is special as Bayes predictive densities are necessarily {\em plug-in} predictive densities, as expanded upon for exponential families in the Appendix.  This represents a slight extension of a result due to Yanigimoto and Ohnishi (2009).  Moreover, the duality with squared error loss opens the way for Bayesian dominance results.  For Kullback-Leibler loss, two  dominance findings are obtained in Subsection 4.2.  For $p=1$, making use of Section 2's representations, we show that the Bayes predictive density $\hat{q}_{\pi_{U,A}}$ improves on $\hat{q}_{mre}$ under Kullback-Leibler loss for both $\theta_1 \geq \theta_2$ or $|\theta_1 - \theta_2| \leq m$.  For the former case, the dominance result is further proven in Theorem \ref{persistence} to be robust with respect to various misspecifications of $\sigma^2_1, \sigma^2_2, $ and $\sigma^2_Y$. Finally, numerical illustrations are presented and commented upon in Section 5. 
 

\section{Bayesian predictive density estimators and skewed normal type distributions}

\subsection{Bayesian predictive density estimators}

\noindent  We provide here a general representation of the Bayes predictive density estimator
of the density of $Y_1$ in model (\ref{model}) associated with a uniform prior on the additional information set $A$.  Multivariate normal priors truncated to $A$ are plausible choices that are also conjugate, lead to similar results, but will not be further considered in this manuscript.  Throughout this manuscript, starting with the next result, we denote $\phi$ as the $N_p(0,I_p)$ p.d.f. 

\begin{lemma}
\label{bayesgeneral}
Consider model (\ref{model}), a Bayes predictive density $\hat{q}_{\pi}$ with respect to prior $\pi$ for $\theta$, and the Bayes predictive density   
$\hat{q}_{\pi_{U,A}}$ with respect to the (uniform) prior $\pi_{U,A}(\theta)= \mathbb{I}_A(\theta_1-\theta_2)$ for $\alpha$-divergence loss $L_{\alpha}$ in (\ref{alphadivergence}).  

\begin{enumerate}
\item[ {\bf (a)}]  For $-1 \leq \alpha < 1$, we have
\begin{equation}
\label{bayesalphageneral}
\hat{q}_{\pi_{U,A}}(y_1;x) \propto \hat{q}_{\hbox{mre}}(y_1;x_1) \, I^{\frac{2}{1-\alpha}}(y_1;x)\,,
\end{equation}
with $\hat{q}_{\hbox{mre}}(y_1;x_1)$ the minimum risk predictive density estimator based on $x_1$ given by a $N_p(x_1, (\sigma^2_1 \frac{(1-\alpha)}{2} + \sigma^2_Y) I_p)  $ density, and  $I(y_1;x) = \mathbb{P}(T \in A)$, with $T \sim N_p\left( \mu_T, \sigma^2_T I_p\right)$,  $\mu_T=\beta (y_1 - x_1) + (x_1-x_2)$, $\sigma^2_T= \frac{2 \sigma^2_1 \sigma^{2}_Y}{(1-\alpha)\sigma^{2}_1 + 2 \sigma^{2}_Y} + \sigma^2_2 $, and 
$\beta= \frac{(1-\alpha) \sigma^{2}_1}{(1-\alpha)\sigma^{2}_1 + 2 \sigma^{2}_Y }$.

\item[ {\bf (b)}]  For $\alpha=1$ (i.e., reverse Kullback-Leibler loss), we have 
\begin{equation}
\label{bayesrkl}
\hat{q}_{\pi}(y_1;x) \sim \,N_p(\mathbb{E}(\theta_1|x), \sigma^2_Y I_p)\,,
\end{equation}
where $\mathbb{E}(\theta_1|x)$ is the posterior expectation of $\theta_1$. 
\end{enumerate}
\end{lemma}
\noindent {\bf Proof.}  {\bf (a)}  As shown by Corcuera and Giummol\`e (1999), the Bayes predictive density estimator of the density of $Y_1$ in  (\ref{model}) under loss $L_{\alpha}$, $\alpha \neq 1$, is given by
\begin{equation}
\nonumber
\hat{q}_{\pi_{U,A}}(y_1;x) \propto  \left\lbrace  \int_{\mathbb{R}^p} \int_{\mathbb{R}^p}  \phi^{(1-\alpha)/2}(\frac{y_1-\theta_1}{\sigma_Y})   \,  \pi(\theta_1, \theta_2|x) \, d\theta_1  \, d\theta_2 \right\rbrace^{2/1-\alpha}.
\end{equation}
\noindent With prior measure $\pi_{U,A}(\theta)= \mathbb{I}_A(\theta_1-\theta_2)$, we obtain
\begin{equation}
\nonumber
\hat{q}_{\pi_{U,A}}(y_1;x) \propto
\left\lbrace  \int_{\mathbb{R}^p} \int_{\mathbb{R}^p}  \phi(\frac{y_1-\theta_1}{\sqrt{\frac{2}{1-\alpha}\sigma^2_Y}})   \,  \phi(\frac{\theta_1-x_1}{\sigma_1} )\, \phi(\frac{\theta_2-x_2}{\sigma_2} ) \, \mathbb{I}_A(\theta_1-\theta_2) \, \, d\theta_1  \, d\theta_2  \, \right\rbrace^{2/1-\alpha},
\end{equation}
given that $\phi^m(z) \propto \phi(m^{1/2}z)$ for $m >0$. By the decomposition 
\begin{equation}
\nonumber
\frac{\|\theta_1 - y_1 \|^2}{a} + \frac{\|\theta_1 - x_1 \|^2}{b} \,=\, \frac{\|y_1-x_1\|^2}{a+b} \, + \, \frac{\|\theta_1 - w \|^2}{\sigma^2_w}\,,
\end{equation}
with $a=\frac{2\sigma^2_Y}{1-\alpha}$, $b=\sigma^2_1$, and $w= \frac{b y_1 + a x_1}{a+b}= \beta y_1 + (1-\beta) x_1$, $\sigma^2_w= \frac{ab}{a+b} = \frac{2 \sigma^2_1 \sigma^{2}_Y}{2 \sigma^2_Y + (1-\alpha)\sigma^{2}_1} $, we obtain

\begin{eqnarray*}
\hat{q}_{\pi_{U,A}}(y_1;x) &\propto & \phi^{2/(1-\alpha)}(\frac{y_1-x_1}{\sqrt{\frac{2 \sigma^2_Y}{1-\alpha} + \sigma^2_1}}) \, \left\lbrace\int_{\mathbb{R}^{2p}} \, \phi(\frac{\theta_1-w}{\sigma_w} ) \, \phi(\frac{\theta_2-x_2}{\sigma_2} ) \, \mathbb{I}_A(\theta_1-\theta_2) \, \, d\theta_1  \, d\theta_2 \right\rbrace^{2/1-\alpha}\, \\
\, & \propto &  \hat{q}_{\hbox{mre}}(y_1;x_1) \, \left\lbrace \mathbb{P}(Z_1-Z_2 \in A) 
\right\rbrace^{2/1-\alpha}\,,
\end{eqnarray*}
with $Z_1, Z_2$ independently distributed as $Z_1 \sim N_p(w, \sigma^2_w)$, $Z_2 \sim N_p(x_2, \sigma^2_2)$.  The result follows by setting $T =^d Z_1-Z_2$.

{\bf (b)}  This part is a consequence of Theorem \ref{Bayes=plugin}, which is a general result for exponential families; presented in the Appendix; and which establishes that Bayes predictive densities are necessarily \emph{ plug--in} predictive densities.  See Example \ref{exampleappendix} for details.
\qed

\noindent  The general form of the Bayes predictive density estimator $\hat{q}_{\pi_{U,A}}$ is thus a weighted version of $\hat{q}_{\hbox{mre}}$, with  the weight a multivariate normal probability
raised to the $2/(1-\alpha)^{th}$ power which is a function of $y_1$ and which depends on $x, \alpha, A$.  Observe that the representation applies in the trivial case $A=\mathbb{R}^p$, yielding $ I=1$ and $\hat{q}_{\hbox{mre}}$ as the Bayes estimator.
As expanded on in Subsection \ref{examples}, the densities $\hat{q}_{\pi_{U,A}}$ for Kullback-Leibler loss relate to skew-normal distributions, and more generally to skewed distributions arising from selection (see for instance Arnold and Beaver, 2002;   Arellano-Valle, Branco and Genton, 2006; among others).  Moreover, it is known (e.g. Liseo and Loperfido, 2003) that  posterior distributions present here also relate to such skew-normal type distributions.    
Lemma \ref{bayesgeneral} does not address the evaluation of the normalization constant for the Bayes predictive density $\hat{q}_{\pi_{U,A}}$, but we now proceed with this for the particular cases of Kullback-Leibler and Hellinger losses, and more generally for cases where $\frac{2}{1-\alpha}$ is a positive integer, i.e., $\alpha=1-\frac{2}{n}$ where $n=1,2,\ldots$. 
In what follows, we denote $1_m$ as the $m$ dimensional column vector with components equal to $1$, and $\otimes$ as the usual Kronecker product.

\begin{lemma}
\label{lemmaalpha-n}
For model (\ref{model}), $\alpha-$divergence loss with $n = \frac{2}{1-\alpha} \in \{1,2,\ldots\}$, the Bayes predictive density $\hat{q}_{\pi_{U,A}}(y_1;x)\,, y_1 \in \mathbb{R}^p,$
with respect to the (uniform) prior $\pi_{U,A}(\theta)= \mathbb{I}_A(\theta_1-\theta_2)$, is given by
\begin{equation}
\label{equationalpha-n}
\hat{q}_{\pi_{U,A}}(y_1;x)\, =  \hat{q}_{\hbox{mre}}(y_1;x_1) \, 
\frac{\{\mathbb{P}(T \in A) \}^n}{\mathbb{P}(\cap_{i=1}^n \{Z_i \in A\})}\,,
\end{equation}
\end{lemma}
with $\hat{q}_{\hbox{mre}}(y_1;x_1)$ a $N_p(x_1, (\sigma^2_1/n + \sigma^2_Y)\hbox{I}_p)$ density, $T \sim N_p(\mu_T, \sigma^2_T I_p)$ with $\mu_T= \beta (y_1-x_1) + (x_1-x_2)$, $\sigma^2_T= \sigma^2_2 + n \sigma^2_Y \beta$, $\beta= \frac{\sigma^2_1}{\sigma^2_1 + n \sigma^2_Y} $, and
$Z=(Z_1, \ldots, Z_n)' \sim N_{np}(\mu_Z, \Sigma_Z)$ with $\mu_Z=1_n \otimes (x_1-x_2)$ and $\Sigma_Z=(\sigma_T^2 + \sigma^2_Y \beta^2) I_{np} + (\frac{\beta^2 \sigma^2_1}{n} 1_n 1_n' \otimes \hbox{I}_p )\, $

\begin{remark}
\label{n=1}
The Kullback-Leibler case corresponds to $n=1$ and the above form of the Bayes predictive density simplifies to
\begin{equation}
\label{klgeneralform}
\hat{q}_{\pi_{U,A}}(y_1;x)\, =  \hat{q}_{\hbox{mre}}(y_1;x_1) \, \frac{\mathbb{P}(T \in A)}{\mathbb{P}(Z_1 \in A)}\,,
\end{equation}
\end{remark} 
with $\hat{q}_{\hbox{mre}}(y_1;x_1)$ a $N_p(x_1, (\sigma^2_1 + \sigma^2_Y)\hbox{I}_p)$ density, $T \sim N_p(\mu_T, \sigma^2_T I_p)$ with $\mu_T= \frac{\sigma^2_1}{\sigma^2_1 + \sigma^2_Y}(y_1-x_1) + (x_1-x_2)$ and $\sigma^2_T= \frac{\sigma^2_1\sigma^2_Y}{\sigma^2_1 + \sigma^2_Y} + \sigma^2_2$, and $Z_1 \sim N_p(x_1-x_2, (\sigma^2_1 + \sigma^2_2)I_p)$.
In the univariate case (i.e., $p=1$), $T$ is univariate normally distributed and the expectation and covariance matrix of $Z$ simplify to 
$1_n (x_1-x_2)$ and $ (\sigma_T^2 + \sigma^2_Y \beta^2) I_n \,+\, \beta^2
\frac{\sigma^2_1}{n} \, 1_n 1_n'$ respectively.  Finally, we point out that the diagonal elements of $\Sigma_Z$ simplify to $\sigma_1^2 + \sigma_2^2$, a result which will arise below several times.

\noindent  {\bf Proof of Lemma \ref{lemmaalpha-n}.}  It suffices to evaluate the normalization constant (say $C$) for the predictive density in (\ref{bayesalphageneral}).  We have   
\begin{eqnarray*}
\label{tobeproven}
C &=& \int_{\mathbb{R}^p} \hat{q}_{\hbox{mre}}(y_1;x_1) \, \{\mathbb{P}(T \in A)\}^n \, dy_1 \\
\, & = & \int_{\mathbb{R}^p} \hat{q}_{\hbox{mre}}(y_1;x_1) \, \mathbb{P} \left(\cap_{i=1}^n \{T_{i} \in A\} \right)\, dy_1\,,
\end{eqnarray*}
with $T_{1}, \ldots ,T_{n}$ independent copies of $T$.  With the change of variables $u_0= \frac{y_1-x_1}{\sqrt{\sigma^2_1/n + \sigma^2_Y}}$ and letting $U_0, U_1, \ldots, U_n$ i.i.d. $N_p(0,I_p)$, we obtain
\begin{eqnarray*}
C \, &=& \, \int_{\mathbb{R}^p} \, \phi(u_0) \,  \mathbb{P}\left( \cap_{i=1}^n \{\sigma_T U_i + \beta u_0 \sqrt{\sigma^2_1/n + \sigma^2_Y} + x_1 - x_2\} \in A \right)\, du_0 \\
\, & = &  \mathbb{P}\left( \cap_{i=1}^n \{\sigma_T U_i + \beta U_0 \sqrt{\sigma^2_1/n + \sigma^2_Y} + x_1 - x_2\} \in A \right)\,, \\
\, & = &  \mathbb{P} \left(\cap_{i=1}^n \{Z_i \in A \} \right).
\end{eqnarray*}
The result follows by verifying that the expectation and covariance matrix of $Z=(Z_1, \ldots, Z_n)'$ are as stated. \qed

\noindent  
The next result presents a useful posterior distribution decomposition, with an accompanying representation of the posterior expectation $\mathbb{E}(\theta_1|x)$ in terms of a truncated multivariate normal expectation.   The latter characterizes the Bayes predictive density under reverse Kullback-Leibler loss in accordance with Lemma \ref{bayesgeneral}, as well as coincide with the expectation under the Bayes Kullback-Leibler predictive density $\hat{q}_{\pi_{U,A}}$.
Specific examples will be presented in Subsection 2.3.4.

\begin{lemma}
\label{posterior}
Consider $X|\theta$ as in model (\ref{model}) and the uniform prior $\pi_{U,A}(\theta)= \mathbb{I}_{A}(\theta_1 - \theta_2)$.  Set  $r= \frac{ \sigma_2^2 }{\sigma_1^2}$, $\omega_1=\theta_1-\theta_2$, and $\omega_2=r\theta_1+\theta_2$.  Then, conditional on $X=x$, $\omega_1$ and $\omega_2$ are independently distributed with
\begin{equation}
\nonumber
\omega_1 \sim N_p(\mu_{\omega_1}, \tau_{\omega_1}^2)  \hbox{ truncated to } A, \;\;\;\; \;\;\;\; \;\;\;\; \omega_2 \sim N_p(\mu_{\omega_2}, \tau_{\omega_2}^2)\,,
\end{equation}
$\mu_{\omega_1}= x_1-x_2$, 
 $\mu_{\omega_2}= rx_1+x_2$, 
 $\tau_{\omega_1}^2= \sigma_1^2 + \sigma_2^2 $, and $\tau_{\omega_2}^2 = 2  \sigma_2^2. $  Furthermore, we have $\mathbb{E}(\theta_1|x) = \frac{1}{1+r} \left(  \mathbb{E}(\omega_1|x) + \mu_{\omega_2}\right)$.
 
\end{lemma}
\noindent {\bf Proof.}   
With the posterior density $\pi(\theta|x) \propto \phi(\frac{\theta_1 - x_1}{\sigma_1})  \; \phi(\frac{\theta_2 - x_2 }{\sigma_2}) \,\, \mathbb{I}_{A}(\theta_1 - \theta_2)$,
the result follows by transforming to $(\omega_1, \omega_2)$.
\qed

\subsection{Examples of Bayesian predictive density estimators}
\label{examples}
\noindent 
With the presentation of the Bayes predictive estimator $\hat{q}_{\pi_{U,A}}$ in Lemmas \ref{bayesgeneral} and \ref{lemmaalpha-n}, which is quite general with respect to the dimension $p$, the additional information set $A$, and the $\alpha-$divergence loss, it is pertinent and instructive to continue with some illustrations.  Moreover, various skewed-normal or skewed-normal type, including new extensions, arise as predictive density estimators.  Such distributions have indeed generated much interest for the last thirty years or so, and continue to do so, as witnessed by the large literature devoted to their study. The most familiar choices of $\alpha-$divergence loss are Kullback-Leibler and Hellinger (i.e., $n=\frac{2}{1-\alpha}=1,2$ below) but the form of the Bayes predictive density estimator $\hat{q}_{\pi_{U,A}}$ is nevertheless expanded upon below in the context of Lemma \ref{lemmaalpha-n}, in view of the connections with an extended family of skewed-normal type distributions (e.g., Definition \ref{skewbala}), which is also of independent interest.  Subsections 2.2.1, 2.2.2, 2.2.3. deal with Kullback-Leibler and $\alpha-$divergence losses for situations: (i) $p=1, A=\mathbb{R}_+$; (ii) $p=1, A=[-m,m]$; (iii) $p \geq 1$ and $A$ a ball of radius $m$ centered at the origin, while Subsection 2.2.4. deals with reverse Kullback-Leibler loss.

\subsubsection{\sc Univariate case with $\theta_1 \geq \theta_2$.}

From (\ref{equationalpha-n}), we obtain for $p=1, A=\mathbb{R}_+$:  $\mathbb{P}(T \in A)= \Phi(\frac{\mu_T}{\sigma_T})$ and 

\begin{equation}
\label{densities}
\hat{q}_{\pi_{U,A}}(y_1;x) \propto  
\frac{1}{\sqrt{\sigma^2_1/n + \sigma^2_Y}} \; \phi(\frac{y_1-x_1}{\sqrt{\sigma^2_1/n + \sigma^2_Y}}) \;\Phi^n(\frac{\beta(y_1-x_1) + (x_1-x_2)}{\sigma_T})\,,
\end{equation}
with $\beta$ and $\sigma^{2}_T$ given in Lemma \ref{lemmaalpha-n}. 
These densities match the following family of densities.

\begin{definition} 
\label{skewbala}
A generalized Balakrishnan type skewed-normal distribution, with shape parameters $n \in \mathbb{N}_+, \alpha_0, \alpha_1 \in \mathbb{R}$, location and scale parameters $\xi$ and $\tau$,  denoted $\mathrm{SN}(n,\alpha_0, \alpha_1, \xi, \tau)$, has density on  $ \mathbb{R}$ given by 
\begin{align}\label{twoSND}
\frac{1}{K_n(\alpha_0, \alpha_1)}\,\frac{1}{\tau}\,\phi(\frac{t-\xi}{\tau})\,\Phi^n(\alpha_0 +\alpha_1 \frac{t-\xi}{\tau})\,,
\end{align}
with 
\begin{align}\label{CDFBivariate}
K_n(\alpha_0, \alpha_1)=\Phi_n\left(\frac{\alpha_0}{\sqrt{1+\alpha_1^2}}, \cdots, \frac{\alpha_0}{\sqrt{1+\alpha_1^2}}; \rho=\frac{\alpha_1^2}{1+\alpha_1^2}\right)\,,
\end{align}
$\Phi_n(\cdot; \rho)$ representing the cdf of a $\mathrm{N}_n(0, \Lambda)$ distribution with covariance matrix $ \Lambda=(1-\rho)\,I_n+ \rho\,{1}_n {1}'_n$.
\end{definition}
\begin{remark}  (The case $n=1$)  \\
$\mathrm{SN}(1,\alpha_0, \alpha_1, \xi, \tau)$ densities are given by  (\ref{twoSND}) with $n=1$ and $K_1(\alpha_0, \alpha_1) = 
\Phi(\frac{\alpha_0}{\sqrt{1+\alpha_1^2}})$.
Properties of $\mathrm{SN}(1,\alpha_0, \alpha_1, \xi, \tau)$ distributions were described by Arnold et al. (1993), as well as Arnold and Beaver (2002), with the particular case $\alpha_0=0$ reducing to the original skew normal density, modulo a location-scale transformation, as presented in Azzalini's seminal 1985 paper.  Namely, the expectation of $T \sim \mathrm{SN}(1, \alpha_0, \alpha_1, \xi, \tau)$ is given by

\begin{equation}
\label{ET}
\mathbb{E}(T) \,=\,  \xi + \tau \frac{\alpha_1}{\sqrt{1+\alpha_1^2}} \, R(\frac{\alpha_0}{\sqrt{1+\alpha_1^2}})\,,
\end{equation}
with $R=: \frac{\phi}{\Phi}$ known as the inverse Mill's ratio.
\end{remark}

\begin{remark}
\label{remarkbala}
For $\alpha_0=0, n=2, 3, \ldots$, the densities were proposed by Balakrishnan as a discussant of Arnold and Beaver (2002), and further analyzed by Gupta and Gupta (2004).    We are not aware of an explicit treatment of such distributions in the general case, but standard techniques may be used to derive the following properties.  For instance, as handled more generally above in the proof of Lemma \ref{lemmaalpha-n},  the normalization constant $K_n$ may be expressed in terms of a multivariate normal c.d.f. by observing that 
\begin{align}
\label{nor-con}
K_n(\alpha_0, \alpha_1)&= \int_{\mathbb{R}} \phi(z) \Phi^n(\alpha_0 + \alpha_1 z ) \,dz \nonumber \\
&=  \mathbb{P}(\cap_{i=1}^n \{U_i \leq \alpha_0  + \alpha_1 U_0  \}) \nonumber \\
&=\mathbb{P}(\cap_{i=1}^n \{W_i \leq \frac{\alpha_0}{\sqrt{1+\alpha_1^2}}  \})\,,
\end{align}
with $(U_0, \ldots, U_n) \sim N_{n+1}(0, I_{n+1})$,  $W_i \stackrel{d}{=}\frac{U_i-\alpha_1\, U_0}{\sqrt{1+\alpha_1^2}}$, for $i=1, \ldots, n$,
 and $(W_1, \ldots, W_n) \sim \mathrm{N}_n (0, \Lambda)$. 
 
In terms of expectation, we have, for $T \sim SN(n,\alpha_0, \alpha_1, \xi, \tau)$, $\mathbb{E}(T)=\xi + \tau \mathbb{E}(W)$ where  $W \sim SN(n,\alpha_0, \alpha_1, 0, 1)$ and
\begin{equation}
\label{ETSN_n}
\mathbb{E}(W) \,=\, \frac{n\alpha_1}{\sqrt{1+\alpha_1^2}} \, \phi(\frac{\alpha_0}{\sqrt{1+\alpha_1^2}})\, \frac{K_{n-1}(\frac{\alpha_0}{\sqrt{1+\alpha_1^2}}, \frac{\alpha_1}{\sqrt{1+\alpha_1^2}})}
{K_n(\alpha_0, \alpha_1)}\,. 
\end{equation}
\end{remark}
This can be obtained via Stein's identity
$\mathbb{E} \, U g(U)=\mathbb{E} g'(U)$ for  differentiable $g$ and $U \sim N(0, 1)$.  Indeed, we have 
\begin{eqnarray*}
\int_{\mathbb{R}} u \phi(u)\, \Phi^n(\alpha_0+\alpha_1 u)\, du  
&=& n \alpha_1 \int_{\mathbb{R}} \phi(u) \phi(\alpha_0+\alpha_1 u)\, \Phi^{n-1}(\alpha_0+\alpha_1 u)\, du \,,
\end{eqnarray*}
and the result follows by making use of the identity  $\phi(u) \, \phi(\alpha_0+\alpha_1 u) \,=\, 
\phi(\frac{\alpha_0}{\sqrt{1+\alpha_1^2}}) \, \phi (v)$, with $v=\sqrt{1+\alpha_1^2} \, u + \frac{\alpha_0 \alpha_1}{\sqrt{1+\alpha_1^2}}$, the change of variables $u \to v$, and the definition of $K_{n-1}$.

\noindent  The connection between the densities of Definition \ref{skewbala} and the predictive densities in (\ref{densities}) is thus explicitly stated as follows, with the Kullback-Leibler and Hellinger cases corresponding to $n=1,2$ respectively.

\begin{corollary}
\label{cor2.1}
For $p=1, A=\mathbb{R}_+$, $\pi_{U,A}(\theta)= \mathbb{I}_A(\theta_1-\theta_2)$, the Bayes predictive density estimator $\hat{q}_{\pi_{U,A}}$ under $\alpha-$divergence loss, with $n=\frac{2}{1-\alpha} \in \mathbb{N_+}$ positive integer,  is given by a $\mathrm{SN}(n,\alpha_0=\frac{x_1-x_2}{\sigma_T}, \alpha_1= \frac{\beta \tau}{\sigma_T}, \xi=x_1, \tau= \sqrt{\frac{\sigma^2_1}{n} + \sigma^2_Y}   )$ density,  with  $\sigma^2_T = \sigma^2_2 + n \beta \sigma_Y^2$ and $\beta = \frac{\sigma^2_1}{\sigma^2_1 + n\sigma^2_Y}$.
\end{corollary}

\begin{remark} For the equal variances case with $\sigma^2_1=\sigma^2_2=\sigma^2_Y=\sigma^2$, the above predictive density estimator is a $\mathrm{SN}(n,\alpha_0= \sqrt{\frac{n+1}{(2n+1) \sigma}} (x_1-x_2), \alpha_1=\sqrt{\frac{1}{n (2n+1)}}, \xi=x_1, \tau = \sqrt{\frac{n+1}{n}} \sigma  )$ density.
\end{remark}

\subsubsection{\sc Univariate case with $|\theta_1 - \theta_2| \leq m$}

From (\ref{equationalpha-n}),  we obtain for $p=1, A=[-m,m]$:  $\mathbb{P}(T \in A)=\Phi(\frac{\mu_T + m}{\sigma_T}) - \Phi(\frac{\mu_T - m}{\sigma_T})$, and we may write

\begin{equation}
\label{with-m}
\hat{q}_{\pi_{U,A}}(y_1;x) =   \frac{1}{\tau}\phi(\frac{t-\xi}{\tau}) \,\, 
 \frac{\lbrace\Phi(\alpha_0+\alpha_1 \,\frac{t-\xi}{\tau}) - \Phi(\alpha_2+\alpha_1 \,\frac{t-\xi}{\tau})\rbrace ^n}{J_n(\alpha_0, \alpha_1, \alpha_2) }\,,
\end{equation}
with $\xi=x_1, \tau=\sqrt{\sigma_1^2/n + \sigma_Y^2}, \alpha_0= \frac{x_1-x_2+m}{\sigma_T}$, $\alpha_1= \frac{\beta \tau}{\sigma_T}$ $\alpha_2= \frac{x_1-x_2-m}{\sigma_T}$, $\beta, \mu_T$, and $\sigma^2_T$ given in Lemma \ref{lemmaalpha-n}, and $J_n(\alpha_0, \alpha_1, \alpha_2)$ (independent of $\xi, \tau$) a special case of the normalization constant given in (\ref{equationalpha-n}).

\noindent   For fixed $n$, the densities in (\ref{with-m}) form a five-parameter family of densities with location and scale parameters $\xi \in \mathbb{R}$ and $\tau \in \mathbb{R}_+$, and shape parameters $\alpha_0, \alpha_1, \alpha_2 \in \mathbb{R}$ such that $\alpha_0 > \alpha_2$.
The Kullback-Leibler predictive densities ($n=1$) match densities introduced by Arnold et al. (1993) with the normalization constant in (\ref{with-m}) simplifying to:
\begin{equation}
\label{J_1}
J_1(\alpha_0, \alpha_1, \alpha_2) \, = \, \Phi(\frac{\alpha_0}{\sqrt{1+\alpha_1^2}}) - \Phi(\frac{\alpha_2}{\sqrt{1+\alpha_1^2}})= \Phi(\frac{m-(x_1-x_2)}{\sqrt{\sigma_1^2+\sigma_2^2}})-\Phi(\frac{-m-(x_1-x_2)}{\sqrt{\sigma_1^2+\sigma_2^2}}).
\end{equation}

The corresponding expectation is readily obtained as in (\ref{ET}) and equals 
\begin{eqnarray}
\nonumber
\mathbb{E}(T) &=& \xi +  \tau  \frac{\alpha_1}{\sqrt{1+\alpha_1^2}} \, \frac{\phi(\frac{\alpha_0}{\sqrt{1+\alpha_1^2}}) - \phi(\frac{\alpha_2}{\sqrt{1+\alpha_1^2}})}{\Phi(\frac{\alpha_0}{\sqrt{1+\alpha_1^2}}) - \Phi(\frac{\alpha_2}{\sqrt{1+\alpha_1^2}}) } \\
\label{ET2-mm}
& = &  x_1 + \frac{\sigma_1^2}{\sqrt{\sigma_1^2 + \sigma_2^2}} \, 
\frac{\phi(\frac{x_1-x_2+m}{\sqrt{\sigma_1^2 + \sigma_2^2}}) - 
\phi(\frac{x_1-x_2-m}{\sqrt{\sigma_1^2 + \sigma_2^2}})}{\Phi(\frac{x_1-x_2+m}{\sqrt{\sigma_1^2 + \sigma_2^2}}) - 
\Phi(\frac{x_1-x_2-m}{\sqrt{\sigma_1^2 + \sigma_2^2}})} \,,
\end{eqnarray}
by using the above values of $\xi, \tau, \alpha_0, \alpha_1, \alpha_2$.

\noindent Hellinger loss yields the Bayes predictive density in (\ref{with-m})
with $n=2$, and a calculation as in Remark \ref{remarkbala} leads to the evaluation 
\begin{equation}
\nonumber
J_2(\alpha_0, \alpha_1, \alpha_2)\,=\, \Phi_2(\alpha_0',\alpha_0'; \alpha_1') + 
\Phi_2(\alpha_2',\alpha_2'; \alpha_1') - 2 \Phi_2(\alpha_0',\alpha_2'; \alpha_1')
\end{equation}
with $\alpha_i'= \frac{\alpha_i}{\sqrt{1+\alpha_1^2}}$ for $i=0,1,2$.

\subsubsection{\sc Multivariate case with $||\theta_1 - \theta_2|| \leq m$.}

\noindent  For $p \geq 1$, the ball $A=\{t \in \mathbb{R}^p:  ||t|| \leq m \}$,
$\mu_T$ and $\sigma_T^2$ as given in Lemma \ref{with-m}, the Bayes predictive density in (\ref{equationalpha-n}) under $\alpha-$divergence loss with $\frac{2}{1-\alpha}=n \in \mathbb{N}_+$ is expressible as
\begin{equation}
\nonumber
\hat{q}_{\pi_{U,A}} \, \propto \, \hat{q}_{\hbox{mre}}(y_1;x_1) \, \{\mathbb{P}(||T||^2 \leq m^2)\}^n\,
\end{equation}
with $T \sim \sigma_T^2 \chi^2_p (\|\mu_T\|^2/\sigma_T^2)$, i.e., the weight attached to $\hat{q}_{\hbox{mre}}$ is proportional to the $n^{th}$ power
of the c.d.f. of a non-central chi-square distribution.  

\noindent  For Kullback-Leibler loss, we obtain from (\ref{equationalpha-n})
\begin{eqnarray}
\nonumber
\hat{q}_{\pi_{U,A}}(y_1;x)\, & = &  \hat{q}_{\hbox{mre}}(y_1;x_1) \,\frac{\mathbb{P}(||T||^2 \leq m^2)}{\mathbb{P}(||Z_1||^2 \leq m^2)}\, \\
\label{extension}
\, & = & \hat{q}_{\hbox{mre}}(y_1;x_1) \,  \frac{\mathbb{F}_{p, \lambda_1(x,y_1)} (m^2/\sigma_T^2)}{\mathbb{F}_{p, \lambda_2(x)} (m^2/(\sigma_1^2+\sigma_2^2))}\,,
\end{eqnarray}
where $F_{p,\lambda}$ represents the c.d.f. of a $\chi^2_p(\lambda)$ distribution,  $\lambda_1(x,y_1) = \frac{\|\mu_T\|^2}{\sigma_T^2}\,=\, \frac{\|\beta (y_1 - x_1) + (x_1 - x_2)\|^2}{\sigma_T^2}$; with $\beta= \frac{\sigma_1^2}{\sigma_1^2 + \sigma_Y^2}, \sigma_T^2= \sigma_2^2 + \beta \sigma_Y^2$;  and 
$\lambda_2(x)= \frac{\|x_1-x_2\|^2}{\sigma_1^2 + \sigma_2^2}$.  
Observe that the non-centrality parameters $\lambda_1$ and $\lambda_2$  
are random, and themselves non-central chi-square distributed as
$\lambda_1(X,Y_1) \sim \chi^2_p(\frac{||\theta_1-\theta_2||^2}{\sigma_T^2})$ and $\lambda_2(X) \sim \chi^2_p(\frac{||\theta_1-\theta_2||^2}{\sigma_1^2+\sigma_2^2})$.
Of course, the above predictive density (\ref{extension}) matches the Kullback-Leibler predictive density given in (\ref{with-m})  for $n=1$, and represents an otherwise interesting multivariate extension.  

\subsubsection{\sc reverse kullback-leibler loss}

It follows from Lemma \ref{bayesgeneral} and Lemma \ref{posterior} (also see Lemma \ref{bayesC})
that the Bayes predictive density estimator $\hat{q}_{\pi_{U,A}}$ for reverse Kullback-Leibler loss, is given by a $N_p(\mathbb{E}(\theta_1|x), \sigma^2_Y I_p)$ density with 
\begin{equation}
\label{posteriorexample}
\mathbb{E}(\theta_1|x)=  \frac{1}{1+r} \, ( \mathbb{E}(\omega_1|x) +  rx_1 + x_2),\, \hbox{ with } \omega_1 \sim N_p(x_1-x_2, (\sigma^2_1 + \sigma^2_2) I_p) \, \hbox{ truncated to } A\,.
\end{equation}
Truncated normal distributions and their expectations are familiar quantities and thus provide expressions for such predictive densities. Alternatively, as mentioned in the paragraph preceding Lemma \ref{posterior}, the expectation  
$\mathbb{E}(\theta_1|x)$ also matches the expected value under the Kullback-Leibler Bayes predictive density $\hat{q}_{U,A}$.  We illustrate these two above approaches by evaluating (\ref{posteriorexample}) for the following situations.   

\begin{enumerate}
\item[ {\bf (I)}]   Consider  $p=1, A=\mathbb{R}_+$ and let $T \sim \hat{q}_{\pi_{U,A}}$ corresponding to Kullback-Leibler loss.  Then, we have 
\begin{equation}
\nonumber
\mathbb{E}(\theta_1|x)= \mathbb{E}(T) =  x_1 + \frac{\sigma_1^2}{\sqrt{\sigma_1^2 + \sigma_2^2}} \, R(\frac{x_1-x_2}{\sqrt{\sigma_1^2 + \sigma_2^2}})\,,
\end{equation}
by using directly (\ref{ET}) and Corollary \ref{cor2.1}.

\item[ {\bf (II)}]  Similarly, for  $p=1, A=[-m,m]$, letting  let $T \sim \hat{q}_{\pi_{U,A}}$ corresponding to Kullback-Leibler loss, we have  $\mathbb{E}(\theta_1|x)= \mathbb{E}(T)$ as given in (\ref{ET2-mm}).
 
\item[ {\bf (III)}]  Consider the ball $A=\{t \in \mathbb{R}^p: \|t \| \leq m\}$ with $p \geq 1$.  Observe that $\mathbb{E}(\omega_1|x)\,=\, \delta_{\pi_{U,A}}(x')$, with $x'=x_1-x_2$, is the Bayes point estimator under squared error loss based on  the model $X' \sim N_p(\mu, (\sigma_1^2 + \sigma_2^2) I_p)$ and the prior $\pi_{U,A}$.  Such an estimator was expressed in terms of the $\chi^2_p(\lambda)$ c.d.f. $F_{p,\lambda}$ by Marchand and Perron (2001, Remark 1).  From their formula and the above connection, we obtain an evaluation of (\ref{posteriorexample}) with

$$ \mathbb{E}(\omega_1|x) \,=\, (x_1-x_2) \, \frac{F_{p+2,\frac{\|x_1-x_2\|^2}{\sigma_1^2 + \sigma_2^2} }(\frac{m^2}{\sigma_1^2 + \sigma_2^2})}{F_{p, \frac{\|x_1-x_2|^2}{\sigma_1^2 + \sigma_2^2}}(\frac{m^2}{\sigma_1^2 + \sigma_2^2})}\,. $$

\end{enumerate}

\section{General dominance results}

We exploit different channels to obtain predictive density estimation improvements on benchmark procedures such as the maximum likelihood predictive density estimator  $\hat{q}_{\hbox{mle}}$ and the minimum risk equivariant predictive density $\hat{q}_{\hbox{mre}}$.  These predictive density estimators are members of the larger class of densities
\begin{equation}
\label{scaleexpansion2}
q_{\hat{\theta}_1,c} \sim  N_p(\hat{\theta}_1(X), c \sigma_Y^2 I_p)\,,
\end{equation}
with, for instance, the choice $\hat{\theta}_1(X)=\hat{\theta}_{1, \hbox{mle}}(X), c=1$ yielding $\hat{q}_{\hbox{mle}}$, and $\hat{\theta_1}(X)=X, c=1+ \frac{(1-\alpha)  \sigma_1^2}{2\sigma_Y^2}$ yielding $\hat{q}_{\hbox{mre}}$ for loss $L_{\alpha}$.
Two main strategies are exploited to produce improvements:  {\bf (A)}  scale expansion and {\bf (B)}  point estimation duality.  

\begin{enumerate}
\item[ {\bf (A)}]  \emph{ Plug--in} predictive densities $q_{\hat{\theta}_1,1}$ were shown in Fourdrinier et al. (2011), in models where $X_2$ is not observed and for Kullback-Leibler loss, to be universally deficient and improved upon uniformly in terms of risk by a subclass of scale expansion variants $q_{\hat{\theta}_1,c}$ with $c-1$ positive and bounded above by a constant depending on the infimum mean squared error of $\hat{\theta_1}$.
An adaptation of their result leads to dominating predictive densities of 
$\hat{q}_{\hbox{mle}}$, as well as other \emph{plug--in} predictive densities which exploit the additional information $\theta_1 - \theta_2 \in A$,  in terms of Kullback-Leibler risk.     This is expanded upon in Subsection \ref{subsection3.1}.

\item[ {\bf (B)}]  

By duality, we mean that the frequentist risk performance of a predictive density $q_{\hat{\theta}_1,c}$  is equivalent to the point estimation frequentist risk of $\hat{\theta}_1$ in estimating $\theta_1$ under an associated dual loss (e.g., Robert, 1996).  For Kullback-Leibler risk, the dual loss is squared error (Lemma \ref{klrkldual}) and our problem connects to the problem of estimating $\theta_1$ with $\theta_1 - \theta_2 \in A$ based on model (\ref{model}).  In turn, as expanded upon in Marchand and Strawderman (2004), improvements for the latter problem can be generated via the rotation technique (Blumenthal and Cohen, 1968, Cohen and Sackrowitz, 1970, van Eeden and Zidek,  2001, 2003) by improvements for a related restricted parameter space problem.  Details are provided in Subsection \ref{subsection3.2}.  \\

\noindent Similarly, for $\alpha-$divergence loss with $\alpha \in (-1,1)$, the predictive density risk performance of $q_{\hat{\theta}_1,c}$ connects  to the point estimation frequentist risk of   $\hat{\theta}_1$ in estimating $\theta_1$, with $\theta_1 - \theta_2 \in A$ based on model (\ref{model}), 
under reflected normal loss $L_{\gamma_0}$ as seen in Lemma \ref{alphadual} below.
In turn, one can capitalize on a result of Kukobawa, Marchand and Strawderman (2015) which provides a sufficient condition, expressed in terms of a dominance condition under squared error loss, for estimator $\hat{\theta}_{1,A}$ to dominate estimator $\hat{\theta}_{1,B}$ under loss $L_{\gamma_0}$. Then, proceeding as above, this latter problem connects to a restricted parameter space and analysis at this lower level provides results all the way back to the original predictive density estimation problem.  Details and illustrations are provided in Subsection \ref{subsection3.2}. 
\end{enumerate}

\subsection{Improvements by variance expansion}
\label{subsection3.1}
Improvements on \emph{ plug--in} predictive density estimators by variance expansion stem from the following result.  

\begin{lemma}
\label{c}
Consider model (\ref{model}) with $\theta_1-\theta_2 \in A$, a given estimator 
$\hat{\theta}_1$ of $\theta_1$, and the problem of estimating the density of $Y_1$ under Kullback-Leibler loss by a predictive density estimator $q_{\hat{\theta}_1,c}$ as in (\ref{scaleexpansion2}).  Let $\underline{R}= \inf_{\theta} \{\mathbb{E}_{\theta}[\|\hat{\theta}_1(X) - \theta_1 \|^2]\}/(p \sigma_Y^2 )$, where the infimum is taken over the parameter space, i.e. $\{\theta \in \mathbb{R}^{2p}: \theta_1 - \theta_2 \in A\}$, and suppose that $\underline{R}>0$.
\begin{enumerate}
\item[{\bf (a)}]
Then, $q_{\hat{\theta}_1,1}$ is inadmissible and dominated by $q_{\hat{\theta}_1,c}$ for $1 < c < c_0(1+\underline{R})$, with $c_0(s)$, for $s>1$, the root $ c \in (s, \infty)$ of $G_s(c)=(1-1/c) \, s- \log c$.  
\item[{\bf (b)}]
Furthermore, we have $s^2 < c_0(s) < e^s$ for all $s >1$, as well as $\lim_{s \to \infty} c_0(s)/e^s =1$. 
\end{enumerate}
\end{lemma}
{\bf Proof}. See Fourdrinier et al. (2011, Theorem 5.1) for part {\bf (a)}.  
For the first part of {\bf (b)}, it suffices to show that {\bf (i)} $G_s(s^2)>0$ and {\bf (ii)} $G_s(e^s)<0$, given that $G_s(\cdot)$ is, for fixed $s$, a decreasing function on $(s, \infty)$.  We have indeed $G_s(e^s)=-se^{-s}<0$, while $G_s(s^2)|_{s=1}
=0$ and $\frac{\partial}{\partial s} G_s(s^2) = (1 - 1/s)^2>0$, which implies {\bf (i)}. 
Finally, set $k_0(s)= \log c_0(s), s>1,$ and observe that the definition of $c_0$ implies that $u(k_0(s))=\frac{k_0(s)}{1 - e^{-k_0(s)}}=s$.  Since $u(k)$ increases in $k \in (1,\infty)$, it must be the case that $k_0(s)$ increases in $s \in (1,\infty)$ with $\lim_{s \to \infty} k_0(s) \geq \lim_{s \to \infty} \log s^2 = \infty.$  The result thus follows since $\lim_{s \to \infty} k_0(s)/s=\lim_{s \to \infty} (1- e^{-k_0(s)})=1$.

\qed

\begin{remark}
\label{completeclass}
Part (b) above is indicative of the large allowance in the degree of expansion that leads to improvement on the \emph{ plug--in} procedure.  However, among these improvements $c \in (1, c_0(1+\underline{R}))$ on $q_{\hat{\theta}_1,1}$, a complete subclass is given by the choices $c \in [1+\underline{R}, c_0(1+\underline{R}))$, while a minimal complete subclass of predictive density estimators $q_{\hat{\theta}_1,c}$ corresponds to the choices $c \in [1+\underline{R}, 1+ \overline{R} ]$, with $\overline{R}= \sup_{\theta} \{\mathbb{E}_{\theta}[\|\hat{\theta}_1(X) - \theta_1 \|^2]\}/(p \sigma_Y^2 )$, where the supremum is taken over the restricted parameter space, with $\theta_1 - \theta_2 \in A\}$  (see Fourdrinier et al., 2011, Remark 5.1).
\end{remark}

\noindent  The above result is, along with Corollary \ref{cormledominance} below, universal with respect to the choice of the \emph{plug-in} estimator $\hat{\theta}_1$, the dimension $p$ and the constraint set $A$.  
We will otherwise focus below on the \emph{plug-in} maximum likelihood predictive density estimator $\hat{q}_{\hbox{mle}}$.  The next result will be used in both this, and the following, subsections.  The first part presents a decomposition of 
 $\hat{\theta}_{1,mle}$, while the second and third parts relate to a  squared error risk decomposition of estimators given by Marchand and Strawderman (2004).

\begin{lemma}
\label{decomposition}
Consider the problem of estimating $\theta_1$ in model (\ref{model}) with $\theta_1-\theta_2 \in A$ and based on $X$.  Set $r=\sigma^2_2 /\sigma^2_1$, 
$\mu_1=(\theta_1 - \theta_2)/(1+r), \mu_2=(r\theta_1+\theta_2)/(1+r)$, 
$W_1=(X_1 - X_2)/(1+r), W_2=(rX_1+X_2)/(1+r)$, and consider the subclass of estimators of $\theta_1$
\begin{equation}
\label{classC}
\mathrm{C}=\left\{ \delta_{\psi}: \, \delta_{\psi}(W_1, W_2)=W_2+\psi(W_1) \right\}\,.
\end{equation}   Then,
\begin{enumerate}
\item[{\bf (a)}]
The maximum likelihood estimator (mle) of $\theta_1$ is a member of $\mathrm{C}$ with $\psi(W_1)$ the mle of $\mu_1$ based on $W_1 \sim N_p(\mu_1, \sigma^2_1/(1+r) I_p)$ and $(1+r) \mu_1 \in A$;  
\item[{\bf (b)}]  The frequentist risk under squared error loss $\|\delta-\theta_1\|^2$ of an estimator $\delta_{\psi} \in \mathrm{C}$ is equal to 
\begin{equation} 
\label{riskdecomposition}
R(\theta, \delta_{\psi})\,=\,  \mathbb{E}_{\mu_1}[ \|\psi(W_1) - \mu_1 \|^2] + \frac{p \sigma^2_2}{1+r}\,;  (1+r) \mu_1 \in A; 
\end{equation}

\item[{\bf (c)}]  Under squared error loss, the estimator $\delta_{\psi_1}$ dominates  $\delta_{\psi_2}$
iff $\psi_1(W_1)$ dominates $\psi_2(W_1)$ as an estimator of $\mu_1$ under loss $\|\psi - \mu_1\|^2$ and the constraint $(1+r) \mu_1 \in A$.
\end{enumerate}
 
\end{lemma}
{\bf Proof.}  Part {\bf (c)} follows immediately from part {\bf (b)}.  As in Marchand and Strawderman (2004), part {\bf (b)} follows since
\begin{eqnarray*}
R(\theta, \delta_{\psi}) & = & \mathbb{E}_{\theta} \left[\|W_2 + \psi(W_1) - \theta_1  \|^2  \right]   \\
\, & = & \,  \mathbb{E}_{\theta} \left[\| \psi(W_1) - \mu_1  \|^2 \, \right] + \mathbb{E}_{\theta} \left[\|W_2  - \mu_2  \|^2 \, \right] \,, 
\end{eqnarray*}
yielding (\ref{riskdecomposition}) given that $W_1$ and $W_2$ are independently distributed with $W_2 \sim N_p(\mu_2, (\sigma^2_2/(1+r)) I_p)$.  
Similarly, for part {\bf (a)}, we have $\hat{\theta}_{1,mle} = \hat{\mu}_{1,mle} +
\hat{\mu}_{2,mle}$ with $\hat{\mu}_{2,mle}(W_1, W_2)=W_2$ and $\hat{\mu}_{1,mle}(W_1,W_2)$ depending only on $W_1 \sim N_p(\mu_1, (\sigma^2_1/(1+r)) I_p)$ given the independence of $W_1$ and $W_2$.  
\qed

\noindent  Combining Lemmas \ref{c} and \ref{decomposition}, we obtain the following.

\begin{corollary}
\label{cormledominance}
Lemma \ref{c} applies to {\em plug-in} predictive density estimators $q_{\delta_{\psi},1} \sim N_p(\delta_{\psi}, \sigma^2_Y I_p)$ with $\delta_{\psi} \in C$, as defined in (\ref{classC}),  and 
\begin{equation}
\label{infimumR}
\underline{R}\,=\, \frac{1}{\sigma^2_Y} \left( \frac{\sigma^2_1 \sigma^2_2 }{\sigma^2_1+\sigma^2_2}   \, + \, \frac{1}{p} \, \inf_{\mu_1}  \mathbb{E}  [\|\psi(W_1) - \mu_1 \|^2]\,\right)\,.  
\end{equation}
Namely, $q_{\delta_{\psi},c} \sim N_p(\delta_{\psi}, c \sigma^2_Y I_p)$ dominates $q_{\delta_{\psi},1}$ for $1 < c < c_0(1+\underline{R})$.  Moreover, we have 
$c_0(1+\underline{R}) \geq (1+ \underline{R} )^2 \geq (1+ \frac{1}{\sigma^2_Y}  \frac{\sigma^2_1 \sigma^2_2 }{\sigma^2_1+\sigma^2_2})^2\,$.   Finally, the above applies to the maximum likelihood  predictive density estimator 
\begin{equation}
\label{qmlerotation}
\hat{q}_{mle} \sim N_p(\hat{\theta}_{1,mle}, \sigma^2_Y I_p)\,, \hbox{ with } 
\hat{\theta}_{1,mle}(X) = W_2 +  \hat{\mu}_{1,mle}(W_1)\,,
\end{equation}
and 
\begin{equation}
\label{infimumRmle}
\underline{R}\,=\, \frac{1}{\sigma^2_Y} \left( \frac{\sigma^2_1 \sigma^2_2 }{\sigma^2_1+\sigma^2_2}   \, + \, \frac{1}{p} \, \inf_{\mu_1}  \mathbb{E}  [\|\hat{\mu}_{1,mle}(W_1) - \mu_1 \|^2]\,\right)\,,
\end{equation}
where $\hat{\mu}_{1,mle}(W_1)$ the mle of $\mu_1$ based on $W_1 \sim N_p(\mu_1, (\sigma^2_1/(1+r)) I_p)$ and under the restriction $(1+r) \mu_1 \in A$.  
\end{corollary}

\noindent  With the above dominance result quite general, one further issue is the determination of the $\underline{R}$, equivalently $c_0(1+\underline{R})$, or a better lower bound.   
Simulation of the mean squared error in (\ref{infimumR}) is a possibility.  Otherwise, analytically, this seems challenging, but the simple univariate order restriction case leads to the following explicit solution.  

\begin{example} 
\label{example-univariatecase}
 (Univariate case with $\theta_1 \geq \theta_2$) \\
Consider model (\ref{model}) with $p=1$ and $A=[0,\infty)$.  The maximum likelihood predictive density estimator $\hat{q}_{mle}$ is given by (\ref{qmlerotation}) with 
$\hat{\mu}_{1,mle}(W_1) = \max(0,W_1)$.  The mean squared error of $\hat{\theta}_{1,mle}(X)$ may be derived from (\ref{riskdecomposition}) as equal to 
\begin{equation}
\nonumber
R(\theta, \hat{\theta}_{1,mle})\,=\,  \mathbb{E}_{\mu_1}[ \, |\hat{\mu}_{1,mle}(W_1) - \mu_1 |^2] + \frac{\sigma^2_2}{1+r}\;, \mu_1 \geq 0.
\end{equation}
A standard calculation for the mle of a non-negative normal mean based on $W_1 \sim N \left(\mu_1, \sigma^2_{W_1}= \sigma^2_1/(1+r)\right)$ yields the expression
\begin{eqnarray*}
\mathbb{E}_{\mu_1}[ \, |\hat{\mu}_{1,mle}(W_1) - \mu_1 |^2] & = & \mu_1^2 \, \Phi(-\frac{\mu_1}{\sigma_{W_1}}) \, + \,  \int_{0}^{\infty} (w_1-\mu_1)^2 \, \phi(\frac{w_1-\mu_1}{\sigma_{W_1}})  \, \frac{1}{\sigma_{W_1}} dw_1 \\
\,  & = &  \sigma^2_{W_1} \left\lbrace \frac{1}{2} \, + \, \rho^2 \Phi(-\rho) \, + \, \int_{0}^{\rho} t^2 \, \phi(t) \, dt \,   \right\rbrace \,,
\end{eqnarray*}
with the change of variables $t=(w_1-\mu_1)/\sigma_{W_1}$, and by setting $\rho=\mu_1/\sigma_{W_1}$.  Furthermore, it is readily verified that the above risk increases in $\mu_1$; as $\frac{d}{d \rho} \left\lbrace \rho^2 \Phi(-\rho) \, + \, \int_{0}^{\rho} t^2 \, \phi(t) \, dt \,  \right\rbrace = 2 \rho \Phi(-\rho) > 0$ for $\rho>0$, ranging from a minimum value of $\sigma^2_{W_1}/2$ to a supremum value of $\sigma^2_{W_1}$.    

Corollary \ref{cormledominance} thus applies with

$$  \underline{R}\,=\, \frac{1}{\sigma^2_Y} ( \frac{\sigma^2_1 \sigma^2_2}{\sigma^2_1 + \sigma^2_2}+ \frac{\sigma^2_{W_1}}{2})  \,  = \,   \frac{\sigma^2_1}{\sigma^2_Y(\sigma^2_1 + \sigma^2_2)} \, (\sigma^2_2 +\sigma^2_1/2)\,. $$
Similarly, Remark \ref{completeclass} applies with $\overline{R} = \sigma^2_1/\sigma^2_Y$.

As a specific illustration of Corollary \ref{cormledominance} and Remark \ref{completeclass}, consider the equal variances case with $\sigma^2_1=\sigma^2_2= \sigma^2_Y$ for which the above yields $ \underline{R}=3/4,  \overline{R}=1$ and for which we can infer that:

\begin{enumerate}
\item[ {\bf (a)}]   $q_{\hat{\theta}_{1,mle}, c}$  dominates $\hat{q}_{mle}$ under Kullback-Leibler loss for $1 < c < c_0(7/4) \approx 3.48066 $ \;

\item[ {\bf (b)}]   Among the class of improvements in {\bf (a)}, the choices 
$7/4 \leq c < c_0(7/4)$ form a minimal complete subclass;

\item[ {\bf (c)}]  A minimal complete subclass among the $q_{\hat{\theta}_{1,mle}, c}$'s is given by the choices 
$c \in [1+\underline{R}, 1+ \overline{R}]=[7/4, 2]$.

\end{enumerate}
\end{example}

\subsection{Improvements through duality}
\label{subsection3.2}

We consider again here predictive density estimators $q_{\hat{\theta_1},c}$, as in 
(\ref{scaleexpansion2}), but focus rather on the role of the plugged-in estimator $\hat{\theta}_1$.   We seek improvements on benchmark choices such as $\hat{q}_{mre}$, and \emph{plug--in} predictive densities with $c=1$.   We begin with known duality results, and namely  Kullback-Leibler and reverse Kullback-Leibler losses which relate to a dual squared error loss.

\begin{lemma}
\label{klrkldual}
For model (\ref{model}), the frequentist risk of the predictive density estimator 
$q_{\hat{\theta}_1, c}$ of the density of $Y_1$,
under both Kullback-Leibler and reverse Kullback-Leibler losses, is dual to the frequentist risk of $\hat{\theta}_1(X)$ for estimating $\theta_1$ under squared error loss $\|\hat{\theta}_1-\theta_1  \|^2$.
 Namely, $q_{\hat{\theta}_{1,A}, c} $ dominates $q_{\hat{\theta}_{1,B}, c} $ under loss $L_{\alpha}$ iff $\hat{\theta}_{1,A}(X)$  dominates  $\hat{\theta}_{1,B}(X)$ under squared error loss. 
\end{lemma}
{\bf Proof.}  We refer to Fourdrinier et al. (2011) for the case of Kullback-Leibler loss.  For reverse Kullback-Leibler loss, the result follows as an application of Theorem \ref{Bayes=plugin-D}; which is a general result  for exponential families presented in the Appendix, and  expanded upon with Example \ref{exampleappendix}.  \qed 

For other $\alpha-$divergence losses, it is reflected normal loss (defined below) which is dual, as shown by Ghosh, Mergel and Datta (2008) for \emph{plug--in} predictive density estimators, as well as scale expansions in (\ref{scaleexpansion2}).   

\begin{lemma} (Duality between $\alpha-$divergence and reflected normal losses) \\
\label{alphadual}
For model (\ref{model}), the frequentist risk of the predictive density estimator 
$q_{\hat{\theta}_1, c}$ of the density of $Y_1$
under $\alpha-$divergence loss (\ref{alphadivergence}), with $|\alpha|<1$, is dual to the frequentist risk of $\hat{\theta}_1(X)$ for estimating $\theta_1$ under reflected normal loss 
\begin{equation}
\label{rnl}
L_{\gamma_0}(\theta_1, \hat{\theta}_1)= 1 - e^{- \|\hat{\theta}_1-\theta_1  \|^2 /2\gamma_0}\,,
\end{equation}
with $\gamma_0=  (\frac{c}{1+\alpha} + \frac{1}{1-\alpha}) \, \sigma_{Y}^2$.   Namely, $q_{\hat{\theta}_{1,A}, c} $ dominates $q_{\hat{\theta}_{1,B}, c} $ under loss $L_{\alpha}$ iff $\hat{\theta}_{1,A}(X)$  dominates  $\hat{\theta}_{1,B}(X)$ under loss $L_{\gamma_0}$ as above. 
\end{lemma}
{\bf  Proof.}  See for instance Marchand, Perron and Yadegari (2017), or again Ghosh, Mergel and Datta (2008). \qed

\begin{remark}
\label{lim}
Observe that $ \lim_{\gamma_0 \to \infty} 2 \gamma_0 \, L_{\gamma_0} (\theta_1, \hat{\theta}_1) \, =\, \|\hat{\theta}_1 - \theta_1 \|^2$, so that the point estimation performance of $\theta_1$ under reflected normal loss $L_{\gamma_0}$
should be expected to match that of squared error loss when $\gamma_0 \to \infty$.  In view of Lemma \ref{klrkldual} and Lemma \ref{alphadual}, this in turn suggests that the $\alpha-$divergence performance of $\hat{q}_{\hat{\theta}_1, c}$ will match both the Kullback-Leibler and reverse  Kullback-Leibler performance when $|\alpha| \to 1$.  Finally, we point  out that the boundedness nature of the loss in (\ref{rnl}) stands out in contrast to 
cases $\alpha=-1,1$ with the unbounded dual squared-error loss.
\end{remark}

Now, pairing Lemma {\ref{klrkldual} and Lemma \ref{decomposition} leads immediately to the following general dominance result for Kullback-Leibler and reverse Kullback-Leibler losses.

\begin{proposition}
\label{propositionKL-RKL}   Consider model (\ref{model}) with $\theta_1 - \theta_2 \in A$ and the problem of estimating the density of $Y_1$ under either Kullback-Leibler or reverse Kullback-Leibler losses.  Set $r=\sigma^2_2 /\sigma^2_1$, 
$W_1=(X_1 - X_2)/(1+r), W_2=(rX_1+X_2)/(1+r), \mu_1=(\theta_1 - \theta_2)/(1+r)$, 
and further consider the subclass of predictive densities $q_{\delta_{\psi},c}$, as in (\ref{scaleexpansion2}) for fixed $c$, with $\delta_{\psi}$ an estimator of $\theta_1$ of the form $\delta_{\psi}(W_1,W_2)= W_2 + \psi(W_1)$.  Then, $q_{\delta_{\psi_A},c}$ dominates $q_{\delta_{\psi_B},c}$ if and only if $\psi_A$ dominates $\psi_B$ as an estimator of $\mu_1$ under loss $\|\psi- \mu_1 \|^2$, for 
$W_1 \sim N_p(\mu_1, \frac{\sigma^2_1}{1+r} I_p)$ and the parametric restriction $(1+r) \mu_1 \in A$.
\end{proposition}
{\bf Proof.}   The result follows from Lemma {\ref{klrkldual} and Lemma \ref{decomposition}.  \qed

The above result connects three problems, namely:

\begin{enumerate}
\item[ {\bf (I)}] the efficiency of $q_{\delta_{\psi},c}$ under KL or RKL loss as a predictive density for $Y_1$ with the additional information $\theta_1-\theta_2 \in A$;

\item[ {\bf (II)}]  the efficiency of $\delta_{\psi}(X)$ as an estimator of $\theta_1$ under squared error loss $\|\delta_{\psi} - \theta_1 \|^2$ with the additional information $\theta_1-\theta_2 \in A$;

\item[ {\bf (III)}]  the efficiency of $\psi(W_1)$ for $W_1 \sim N_p(\mu_1, \sigma^2_1/(1+r) I_p)$ as an estimator of $\mu_1$ under squared error loss  $\|\psi-\mu_1 \|^2$ with the parametric restriction $(1+r) \mu_1 \in A$.

\end{enumerate}

Previous authors (Blumenthal and Cohen, 1968; Cohen and Sackrowitz, 1970; van Eeden and Zidek, (2001, 2003), for $p=1$; Marchand and Strawderman, 2004, for $p \geq 1$) have exploited the {\bf (II)-(III)} connection (i.e., Lemma \ref{decomposition}) to obtain findings for problem {\bf (II)} based on restricted parameter space findings for {\bf (III)}.  The above Proposition further exploits connections {\bf (I)-(II)} (i.e., Lemma \ref{klrkldual}) to derive findings for predictive density estimation problem {\bf (I)} from restricted parameter space findings for {\bf (III)}.   Consequently, findings for {\bf (III)-(II)} provide findings for our predictive density estimation problem {\bf (I)}, and we refer for Marchand and Strawderman (2004), as well as the references therein, for examples of such results.  An example, which is also illustrative of $\alpha-$divergence results, is provided below at the end of this section. 

\noindent For $\alpha-$divergence losses other than Kullback-Leibler and reverse Kullback-Leibler, the above scheme is not immediately available for the dual reflected normal loss since Lemma \ref{decomposition} is intimately linked to squared error loss.  However,  a slight extension of Lemma 3.3 of Kubokawa, Marchand and Strawderman (2015); exploiting a concave loss technique dating back to Brandwein and Strawderman (1980); permits us to connect (but only in one direction) reflected normal loss to squared error loss, and consequently the efficiency of predictive densities under $\alpha$-divergence loss to point estimation in restricted parameter spaces as in {\bf (III)} above.

\begin{lemma}
\label{concave}
Consider model (\ref{model}) and the problem of estimating $\theta_1$ based on $X$, with $\theta_1 - \theta_2 \in A$ and reflected normal loss as in (\ref{rnl}) with 
$|\alpha| < 1$.    Then $\hat{\theta}_1(X)$ dominates $X_1$ whenever $\hat{\theta}_1(Z)$ dominates $Z_1$ as an estimate of $\theta_1$, under squared error loss $\|\hat{\theta}_1- \theta_1\|^2$, with  $\theta_1 - \theta_2 \in A$, for the model 
\begin{align}
\label{modelZ}
Z=\displaystyle \binom{Z_1}{Z_2} \sim\mathrm{N}_{2p}\left(\displaystyle \theta=\binom{\theta_1}{\theta_2},\,\Sigma_Z=\bigl( \begin{smallmatrix} 
  \sigma_{Z_1}^2 I_{p} & 0\\
  0 & \sigma_2^2 I_{p} 
\end{smallmatrix} \bigr)\right)\,,
\end{align}
with $\sigma_{Z_1}^2 = \frac{\gamma \sigma_1^2}{\gamma+\sigma_1^2}$. 
\end{lemma}
\noindent {\bf Proof.}
Denote the loss $\rho(\|\hat{\theta}_1 - \theta_1\|^2)$ with $\rho(t)=1- e^{-t/2\gamma}$.  Since $\rho$ is concave, we have for all $x=(x_1,x_2)' \in \mathbb{R}^{2p}$:
\begin{equation}
\nonumber \rho (\|\hat{\theta}_1(x) - \theta_1\|^2) - \rho(\|x_1 - \theta_1\|^2)  
\leq   \rho'(\|x_1 - \theta_1\|^2) \, \left(\|\hat{\theta}_1(x) - \theta_1\|^2
 - \|x_1 - \theta_1\|^2)   \right)\,.
 \end{equation}
\noindent   With $\rho'(t)=\frac{1}{2 \gamma}  e^{-t/2\gamma}$, we have for the difference in risks and $Z \sim f_Z$:
\begin{eqnarray*}
\Delta(\theta) &=& R(\theta, \hat{\theta}_1) - R(\theta, X_1) \\
& \leq &  \frac{1}{2\gamma} \frac{1}{(2 \pi \sigma_1 \sigma_2)^p} \int_{\mathbb{R}^{2p}} e^{-\frac{\|x_1-\theta_1\|^2}{2\gamma}}  \left(\|\hat{\theta}_1(x) - \theta_1\|^2 - \|x_1 - \theta_1\|^2)   \right) \, e^{-\frac{\|x_1-\theta_1\|^2}{2\sigma_1^2} - \frac{\|x_2-\theta_2\|^2}{2\sigma_2^2}} \, dx \\
& = &  \frac{1}{2 \gamma} (\frac{\gamma}{\gamma+\sigma_1^2})^{p/2}  \int_{\mathbb{R}^{2p}}  \left(\|\hat{\theta}_1(z) - \theta_1\|^2 - \|z_1 - \theta_1\|^2)   \right)\, f_Z(z) dz \,,
\end{eqnarray*}
establishing the result. \qed

\begin{proposition}
\label{propositionalpha}   Consider model (\ref{model}) with $\theta_1 - \theta_2 \in A$ and the problem of estimating the density of $Y_1$ under either Kullback-Leibler or reverse Kullback-Leibler losses.  Set $r=\sigma^2_2 /\sigma^2_1$, 
$W_1=(X_1 - X_2)/(1+r), W_2=(rX_1+X_2)/(1+r), \mu_1=(\theta_1 - \theta_2)/(1+r)$, 
and further consider the subclass of predictive densities $q_{\delta_{\psi},c}$, as in (\ref{scaleexpansion2}) for fixed $c$, with $\delta_{\psi}$ an estimator of $\theta_1$ of the form $\delta_{\psi}(W_1,W_2)= W_2 + \psi(W_1)$.  Then, $q_{\delta_{\psi_A},c}$ dominates $q_{\delta_{\psi_B},c}$ as long as $\psi_A$ dominates $\psi_B$ as an estimator of $\mu_1$ under loss $\|\psi- \mu_1 \|^2$, for 
$W_1 \sim N_p(\mu_1, \frac{\sigma^2_{Z_1}}{1+r} I_p)$, the parametric restriction $(1+r) \mu_1 \in A$, and $\sigma^2_{Z_1} = \frac{\{(1+\alpha) + c(1-\alpha)\} \, \sigma^2_1}{\{(1+\alpha) + c(1-\alpha)\} \, + \, (1-\alpha^2) \, \sigma^2_1/\sigma^2_Y}$.
\end{proposition}
{\bf Proof.}   The result follows from Lemma {\ref{alphadual} and its dual reflected normal loss $L_{\gamma_0}$, the use of Lemma \ref{concave} applied to $\sigma^2_{Z_1}= \frac{\gamma_0 \sigma^2_1}{\gamma_0+\sigma^2_1}$, and an application of part {\bf (c)} of Lemma \ref{decomposition} to $Z$ as distributed in (\ref{modelZ}).  \qed

\begin{remark} 
\label{remark |alpha|=1}
 Proposition \ref{propositionalpha} holds as stated for $|\alpha|=1$ and is thus a continuation of the sufficiency part of Proposition \ref{propositionKL-RKL}.  As well, the above result provides positive findings as long as $\psi_B$ is inadmissible under squared error loss and dominating estimators $\psi_A$ are available.  Many particular cases follow from the above.  These include: (i) Hellinger loss with $\alpha=0$ and $\sigma^2_{Z_1}$ simplifying to $\left\lbrace(c+1)/(c+1+\sigma^2_1/\sigma^2_Y) \right\rbrace \, \sigma^2_1$; (ii) {\em plug-in} predictive densities with $c=1$; (iii) cases where $q_{\delta_{\psi_B}} \equiv \hat{q}_{mre}$  with the corresponding choice $c=1+ \frac{(1-\alpha) \sigma^2_1}{2 \sigma^2_Y}$  yielding 

$$\sigma^2_{Z_1} = \frac{4 \sigma^2_Y + (1-\alpha)^2 \sigma^2_1}{4 \sigma^2_Y + 
(3 + \alpha) (1 - \alpha) \, \sigma^2_1} \,\,\,  \sigma^2_1\,.$$ 
\end{remark}

The above $\alpha-$ divergence result connects four problems, namely:

\begin{enumerate}
\item[ {\bf (I)}] the efficiency of $q_{\delta_{\psi},c}$ under $\alpha-$divergence loss, $-1 < \alpha < 1$, as a predictive density for $Y_1$ with the additional information $\theta_1-\theta_2 \in A$;

\item[ {\bf (IB)}] the efficiency of $\delta_{\psi}(X)$ as an estimator of $\theta_1$ under reflected normal loss $L_{\gamma_0}$ with $\gamma_0=  (\frac{c}{1+\alpha} + \frac{1}{1-\alpha}) \, \sigma_{Y}^2$ with the additional information $\theta_1-\theta_2 \in A$;

\item[ {\bf (II)}]  the efficiency of $\delta_{\psi}(Z)$, for $Z$ distributed as in (\ref{modelZ}) with $\sigma^2_{Z_1}= (\gamma_0 \sigma^2_1) / (\gamma_0 + \sigma^2_1)$, as an estimator of $\theta_1$ under squared error loss $\|\delta_{\psi} - \theta_1 \|^2$ with the additional information $\theta_1-\theta_2 \in A$;

\item[ {\bf (III)}]  the efficiency of $\psi(W_1)$ for $W_1 \sim N_p(\mu_1, \sigma^2_{Z_1}/(1+r) I_p)$ as an estimator of $\mu_1$ under squared error loss  $\|\psi-\mu_1 \|^2$ with the parametric restriction $(1+r) \mu_1 \in A$.

\end{enumerate}

\begin{example}
Here is an illustration of both Propositions \ref{propositionKL-RKL} and \ref{propositionalpha}.   Consider model (\ref{model}) with $A$ a convex set with a non-empty interior, and $\alpha-$divergence loss ($|\alpha| \leq 1$) for assessing a predictive density for $Y_1$.   Further consider the minimum risk predictive density $\hat{q}_{mre}$ as a benchmark procedure, which is of the form $q_{\delta_{\psi_B}}$ as in Proposition \ref{propositionalpha} with $\delta_{\psi_B} \in C$, $\psi_B(W_1)=W_1$ and $c=c_{mre}=1+ (1-\alpha) \sigma^2_1/(2 \sigma^2_Y)$.  Now consider the Bayes estimator $\psi_U(W_1)$ under squared error loss of $\mu_1$ associated with a uniform prior on the restricted parameter space $(1+r) \mu_1 \in A$, for $W_1 \sim N_p((\mu_1, \frac{\sigma^2_{Z_1}}{1+r} I_p)$ as in Proposition \ref{propositionalpha}.   It follows from Hartigan's theorem (Hartigan, 2003; Marchand and Strawderman, 2004) that $\psi_A(W_1) \equiv \psi_U(W_1)$ dominates 
$\psi_B(W_1)$ under loss $\|\psi-\mu_1\|^2$ and for $(1+r) \mu_1 \in A$.  It thus follows from Proposition \ref{propositionalpha} that the predictive density $q_{\delta_{\psi_B}, c_{mre}} \sim N_p(\delta_{\psi_B}(X), (\frac{1-\alpha}{2} \sigma^2_1 + \sigma^2_Y) I_p) $   dominates $\hat{q}_{mre}$ under $\alpha-$divergence loss with $\delta_{\psi_B}(X)= \frac{rX_1+X_2}{1+r} + \psi_U(\frac{X_1-X_2}{1+r})$.  The dominance result is unified with respect to $\alpha \in [-1,1]$, the dimension $p$, and the set $A$.
\end{example}

We conclude this section with an adaptive two-step strategy, building on both variance expansion and improvements through duality, to optimise on potential Kullback-Leibler improvements on a maximum likelihood estimator predictive density estimator in model (\ref{model}) of the form $\hat{q}_{mle} \sim N_p(\hat{\theta}_{1,mle}, \sigma^2_Y I_p)$, in cases where point estimation improvements on $\hat{\theta}_{1,mle}(X)$ under squared error loss are readily available.

\begin{enumerate}
\item[ {\bf (I)}]  Select an estimator $\hat{\theta}_1^*$ which dominates $\hat{\theta}_{1,mle}$ under squared error loss.  This may be achieved via part {\bf (c)} of Lemma \ref{decomposition} resulting in a dominating estimator of the form $\hat{\theta}_1^*(X)= W_2 + \psi^*(W_1)= (rX_1+X_2)/(1+r) + \psi^*((X_1-X_2)/(1+r))$ 
where $\psi^*(W_1)$ dominates $\hat{\mu}_{1,mle}(W_1)$ as an estimator of $\mu_1$ under squared error loss and the restriction $(1+r) \mu_1 \in A$.

\item[ {\bf (II)}]  Now, with the {\em plug-in} predictive density estimator $q_{\hat{\theta}_{1^*},1}$ dominating $\hat{q}_{mle}$, further improve $q_{\hat{\theta}_{1^*},1}$ by a variance expanded $q_{\hat{\theta}_{1^*},c}$.  Suitable choices of $c$ are prescribed by Corollary \ref{cormledominance} and given by $c_0(1+\underline{R})$, with $\underline{R}$ given in (\ref{infimumR}).  The evaluation of $\underline{R}$ hinges on the infimum risk $\inf_{\mu_1} \mathbb{E}  [\|\psi^*(W_1) - \mu_1 \|^2]$, and such a quantity can be either estimated by simulation, derived in some cases analytically, or safely underestimated by $0$.   

\end{enumerate}

Examples where the above can be applied include the cases: {\bf (i)}  $A=[0,\infty)$ with the use of Shao and Strawderman's (1996) dominating estimators, and {\bf (ii)}  $A$ the ball of radius $m$ centered at the origin with the use of Marchand and Perron's (2001) dominating estimators.  \footnote{Alternatively, one could expand the variance first, and then improve on the {\em plug-in}; such as using a Shao and Strawderman estimator to obtain an improvement on $q_{\hat{\theta}_{1,mle},c}$ in Example \ref{example-univariatecase}; but this may be suboptimal in view of the complete class considerations of Remark \ref{completeclass}.}

\section{Bayesian dominance results}
\label{bayesdominators}

In the previous section, we studied the efficiency of predictive densities as in (\ref{scaleexpansion2}) and elaborated on methods to obtain improvements, whenever possible, for instance on \emph{plug-in} and minimum risk equivariant predictive density estimators.  We focus here on Bayesian improvements, for  reverse Kullback-Leibler and Kullback-Leibler losses, of the benchmark minimum risk equivariant predictive density estimator. For Kullback-Leibler loss, we establish that the uniform Bayes predictive density estimator $\hat{q}_{\pi_{U,A}}$ dominates $\hat{q}_{mre}$ for the univariate cases where $\theta_1-\theta_2$ is either restricted to a compact interval,  lower-bounded or upper-bounded.
Our findings for reverse Kullback-Leibler loss are more wide ranging.  Indeed, we exploit the fact that Bayes predictive density estimators are \emph{plug--in} predictive density estimators,  that the comparison of such procedures is dual to point estimation comparisons under squared error loss, and that we thus can capitalize on existing results for our purposes via Lemma \ref{decomposition}. Such properties are, as expanded upon in the Appendix, quite general for exponential families and reverse Kullback-Leibler loss.

\subsection{Reverse Kullback-Leibler loss}

We begin with an identification of Bayes predictive densities that belong to the class $C$ given in (\ref{classC}), which will permit us to apply Lemma \ref{decomposition} in decomposing the frequentist risk of such procedures.  
This formalizes and extends representation (\ref{posteriorexample}).

\begin{lemma}
\label{bayesC}
Consider model (\ref{model}) and the problem of estimating $\theta_1$ based on $X$  with $\theta_1 - \theta_2 \in A$ and loss $\|\delta-\theta_1  \|^2$.  
Set $r=\sigma^2_2 /\sigma^2_1$, $\mu_1=(\theta_1 - \theta_2)/(1+r), \mu_2=(r\theta_1+\theta_2)/(1+r)$, $W_1=(X_1 - X_2)/(1+r), W_2=(rX_1+X_2)/(1+r)$, and consider prior densities of the form $\pi(\theta) \, = \, \pi_1(\mu_1) \, \mathbb{I}_{A}((1+r) \mu_1)\,\mathbb{I}_{\mathbb{R}^p}(\mu_2)$.  Then, the 
corresponding Bayes estimators $\hat{\theta}_{1,\pi}$ are members of the subclass $C$, as defined in (\ref{classC}), and are given by
\begin{equation}
\label{hattheta1}
\hat{\theta}_{1,\pi}(X)= \psi_{\pi}(W_1) + W_2\,, 
\end{equation}
where $\psi_{\pi}(W_1)$ is the Bayes estimator  based on $W_1 \sim N_p(\mu_1, \frac{\sigma^2_1}{1+r} I_p)$ of $\mu_1$ for loss $\|\psi- \mu_1 \|^2$ and prior $\pi_1(\mu_1) \, \mathbb{I}_{A}((1+r) \mu_1)\,$.
\end{lemma}
{\bf Proof.}   The result follows since the Bayes point estimator of $\theta_1$ is given by
$\mathbb{E}(\theta_1|x)=\mathbb{E}(\mu_1|w_1,w_2) + \mathbb{E}(\mu_2|w_1,w_2) = \mathbb{E}(\mu_1|w_1) + \mathbb{E}(\mu_2|w_2) = \psi_{\pi}(w_1) + w_2$, given the independence of $W_1,W_2$ and the multiplicative aspect of the prior which imply $\mu_1|w_1,w_2 =^d \mu_1|w_1$ and  
$\mu_2|w_1,w_2 =^d \mu_2|w_1$.  \qed

\begin{proposition}
\label{propositionrkl}
Consider model (\ref{model}) with $\theta_1 - \theta_2 \in A$, a prior density of the form $\pi(\theta) \, = \, \pi_1(\mu_1) \, \mathbb{I}_{A}((1+r) \mu_1)\,$, and the corresponding Bayes predictive density $\hat{q}_{\pi}$ for estimating the density of $Y_1$ under reverse Kullback-Leibler loss.  Set $r=\sigma^2_2 /\sigma^2_1$, 
$W_1=(X_1 - X_2)/(1+r), W_2=(rX_1+X_2)/(1+r), \mu_1=(\theta_1 - \theta_2)/(1+r)$, and let $q_{\delta_{\psi_0}}(\cdot;X) \sim N_p(\delta_{\psi_0}(X), \sigma^2_Y I_p) $ be a competing {\em plug-in} predictive density with $\delta_{\psi_0} \in C$ of the form $\delta_{\psi}(W_1,W_2) = \psi_0(W_1) + W_2$.  Then, $\hat{q}_{\pi}(\cdot;X)$
dominates $q_{\delta_{\psi_0}}(\cdot;X)$ if and only if  the Bayes estimator $\psi_{\pi}(W_1)$, with respect to the prior $\pi_1(\mu_1) \, \mathbb{I}_{A}((1+r) \mu_1)\,$, dominates $\psi_0(W_1)$ as an estimator of $\mu_1$ under loss $\|\psi-\mu_1 \|^2$, for $ W_1 \sim N_p(\mu_1, \frac{\sigma^2_{1}}{1+r} I_p)$ and $(1+r) \mu_1 \in A$. 
\end{proposition}
{\bf Proof.}  Part {\bf (b)}  of Lemma \ref{bayesgeneral} and Lemma \ref{bayesC} tell us that $\hat{q}_{\pi}$ is a plug-in predictive density of the form $N_p(\hat{\theta}_{1,\pi}(X), 
\sigma^2_Y I_p)$ with $\hat{\theta}_{1,\pi}(X)$ as in (\ref{hattheta1}).   In turn, Lemma \ref{klrkldual} implies that the reverse Kullback-Leibler risk comparison between $\hat{q}_{\pi}$ and $q_{\delta_{\psi_0}}$ hinges on the mean squared error comparison between $\hat{\theta}_{1,\pi}$ and $\delta_{\psi_0}$ under model (\ref{model}).  Finally, the result follows by making use of Lemma \ref{decomposition}.  \qed

We pursue with applications.

\begin{example}
Consider the context of Proposition \ref{propositionrkl} with $A$ a convex set with a non-empty interior,  the restricted to $A$ uniform prior $\pi_{U,A}(\theta) = \mathbb{I}_{A}(\theta_1-\theta_2)$ and its corresponding Bayes predictive density $\hat{q}_{\pi_{U,A}}$ (see section 2.2.4.),  and the minimum risk predictive density $\hat{q}_{mre}(\cdot; X) \sim N_p(X_1, \sigma^2_Y I_p)$.  It follows from Hartigan's theorem that the Bayes estimator $\psi_U(W_1)$ dominates $\psi_0(W_1)=W_1$ under squared error loss.  Hence, from Proposition \ref{propositionrkl}, it follows that the
Bayes predictive density $\hat{q}_{\pi_{U,A}}$ dominates $\hat{q}_{mre}$ for reverse Kullback-Leibler loss.  The result is general with respect to the choices of $p$ and $A$.

For $p=1$ and $A=[-m,m]$, Kubokawa (2005), as well as Marchand and Payandeh (2011), provide  alternative Bayes estimators $\psi_{\pi_a}(W_1)$ which dominate as well $W_1$ for priors 
$\pi_a$ supported on the set $\mu_1 \in [\frac{-m}{1+r}, \frac{m}{1+r}]$.  In turn, and as above for the uniform prior, it thus follows that the corresponding Bayes predictive densities $\hat{q}_{\pi}(\cdot;X) \sim N(\psi_{\pi_a}(W_1) + W_2, \sigma^2_Y)$ dominate  $\hat{q}_{mre}$
with $\pi(\theta)=\pi_a(\mu_1) \mathbb{I}_{\mathbb{R}}(\mu_2)$. 
\end{example}

\begin{remark}
\label{Stein+}
For $p \geq 3$, $\hat{q}_{\pi_U}$, as well as {\em plug-in} predictive density of the form $q_{\delta_{\psi_0}}(\cdot;X) \sim N_p(\psi_0(W_1) + W_2, \sigma^2_Y I_p)$, are inadmissible and dominated by predictive densities $q_{\delta_{\psi_0, \psi_1}}(\cdot;X) \sim N_p(\psi_0(W_1) + \psi_1(W_2), \sigma^2_Y I_p)$ where $\psi_1(W_2)$ is an estimator of $\mu_2$, for $W_2 \sim N_p(\mu_2, \frac{\sigma^2_2}{1+r} I_p)$, which dominates $W_2$. Stein estimation findings (e.g., Stein, 1981) provide many such dominating estimators, including Bayesian improvements.  For instance, for $p \geq 3$ and a superharmonic prior $\pi_2$ for $\mu_2$, the predictive density $\hat{q}_{\pi_U}$ is dominated by the Bayes predictive density  $q_{\delta_{\psi_U, \psi_{\pi_2}}}(\cdot;X) \sim N_p(\psi_U(W_1) + \psi_{\pi_2}(W_2), \sigma^2_Y I_p)$, associated with the prior $\pi(\theta)= \, \mathbb{I}_A((1+r)\mu_1)\, \pi_2(\mu_2)$.
The above inferences come about a rewriting of Lemma \ref{decomposition} for estimators of the form $\psi_0(W_1) + \psi_1(W_2)$, with $\psi_0 \equiv \psi_U$ for the case of $\hat{q}_{mre}$ and its use as in Proposition \ref{propositionrkl}.
\end{remark}

\begin{example}
\label{rklmle}
Consider the context of Proposition \ref{propositionrkl} and the maximum likelihood predictive density estimator $\hat{q}_{mle} \sim N_p(\hat{\theta}_{1,mle}, \sigma^2_Y I_p)$  with $\hat{\theta}_{1,mle}(X) = W_2 + \psi_0(W_1)$, as in (\ref{qmlerotation}) with  $\psi_0(W_1)=\hat{\mu}_{1,mle}(W_1)$.  It follows from Lemma \ref{klrkldual} that {\em plug-in} predictive densities $N_p(\psi_1(W_1) + W_2, \sigma^2_Y I_p)$ dominate $\hat{q}_{mle}$ under reverse Kullback-Leibler loss if and only if $\psi_1(W_1)$ dominates $\hat{\mu}_{1,mle}(W_1)$ under squared error loss.  In particular and in accordance with Proposition \ref{propositionrkl},  a Bayes predictive density $\hat{q}_{\pi}$, for prior  $\pi(\theta) = \pi_1(\mu_1) \, \mathbb{I}_A((1+r) \mu_1) \mathbb{I}_{\mathbb{R}^p}(\mu_2)$,  dominates  $\hat{q}_{mle}$ if and only if $\psi_{\pi}(W_1)$ dominates $\hat{\mu}_{1,mle}(W_1)$, where $\psi_{\pi}(W_1)$ is the Bayes point estimator of $\mu_1$ for prior $ \pi_1(\mu_1) \, \mathbb{I}_A((1+r) \mu_1)$.  The determination of such dominating Bayesian $\psi_{\pi}$ is challenging though.  For the specific case of $A$ being a $p$-dimensional ball of radius $m$ centered at the origin,  Marchand and Perron (2001), as well as Fourdrinier and Marchand (2010),  provide several applicable Bayesian dominance results.  

\end{example}

\subsection{Kullback-Leibler loss}

In this subsection, we show, for $\theta_1-\theta_2$ either lower bounded, upper bounded, or bounded to an interval, that the uniform Bayes predictive density estimator 
$\hat{q}_{\pi_{U,A}}$ improves uniformly on the minimum risk equivariant predictive density estimator $\hat{q}_{\hbox{mre}}$ under Kullback-Leibler loss.  Without loss of generality (given Remark \ref{transformation}), we consider the restrictions $\theta_1\geq \theta_2$ and $|\theta_1-\theta_2|\leq m$.  We also investigate situations where the variances of model (\ref{model}) are misspecified, but where the dominance persists.  We begin with the lower bounded case.

\begin{theorem} \label{Dominancecon}
Consider model (\ref{model}) with $p=1$ and $A=[0,\infty)$.  For estimating the density of $Y_1$ under Kullback-Leibler loss, the Bayes predictive density $\hat{q}_{\pi_{U,A}}$ dominates the minimum risk equivariant predictive density estimator $\hat{q}_{\hbox{mre}}$.  The Kullback-Leibler risks are equal iff $\theta_1=\theta_2$.  
\end{theorem}
{\bf Proof.}
Making use of Corollary \ref{cor2.1}'s representation of $\hat{q}_{\pi_{U,A}}$, the  
difference in risks is given by
\begin{eqnarray}
\nonumber \Delta(\theta) &= & R_{KL}(\theta, \hat{q}_{\hbox{mre}}) - R_{KL}(\theta, \hat{q}_{\pi_{U,A}}) \\
\nonumber \, & = &   \mathbb{E}^{X,Y_1} \log \left( \frac{\hat{q}_{\pi_{U,A}}(Y_1;X)}{\hat{q}_{\hbox{mre}}(Y_1;X)} \right) \\
\label{diffrisks}\, & = &  \mathbb{E}^{X,Y_1} \log \left( \Phi(\alpha_0 + \alpha_1 \frac{Y_1 -X_1}{\tau}) \right) 
- \mathbb{E}^{X,Y_1} \log \left( \Phi(\frac{\alpha_0}{\sqrt{1+ \alpha^2}}\right)\,,
\end{eqnarray}
with $\alpha_0=\frac{X_1-X_2}{\sigma_T}$, $\alpha_1 = \frac{\beta \tau}{\sigma_T}$, $\tau = \sqrt{\sigma^2_1 + \sigma^2_Y}$, $\beta=\frac{\sigma^2_1}{\sigma^2_1 + \sigma^2_Y}$, and $\sigma^2_T = \sigma^2_2 + \beta \sigma^2_Y$.  
Now, observe that
\begin{equation}
\alpha_0 + \alpha_1 \frac{Y_1 -X_1}{\tau} = \frac{1}{\sigma_T} (X_1-X_2 + \beta (Y_1-X_1) \sim N(\frac{\theta_1-\theta_2}{\sigma_T}, 1)\,,
\end{equation}
and 
\begin{equation}
\frac{\alpha_0}{\sqrt{1+ \alpha_1^2}} = \frac{X_1-X_2}{\sqrt{\sigma_1^2 + \sigma_2^2}} \sim N(\frac{\theta_1-\theta_2}{\sqrt{\sigma_1^2 + \sigma_2^2}}, 1)\,.
\end{equation}
We thus can write  $$\Delta(\theta)\,=\, \mathbb{E} \, G(Z)\,,$$
$$ \hbox{ with }  G(Z)=  \log \Phi(Z+\frac{\theta_1-\theta_2}{\sigma_T}) - \log  \Phi(Z+\frac{\theta_1-\theta_2}{\sqrt{\sigma_1^2 + \sigma_2^2}}) \,,Z \sim N(0,1)\,.$$  
With $\theta_1-\theta_2 \geq 0$ and $\sigma_T^2 < \sigma_1^2 + \sigma_2^2$, we infer that $\mathbb{P}_{\theta}(G(Z)\geq 0)=1$ and $\Delta(\theta)\geq 0$ for all $\theta$ such that $|\theta_1-\theta_2|\leq m$, with equality iff $\theta_1-\theta_2=0$. \qed

We now obtain an analogue dominance result in the univariate case for the additional information $\theta_1 - \theta_2 \in [-m,m]$.  

\begin{theorem} 
\label{dominancedb}
Consider model (\ref{model}) with $p=1$ and $A=[-m,m]$.  For estimating the density of $Y_1$ under Kullback-Leibler loss, the Bayes predictive density $\hat{q}_{\pi_{U,A}}$ (strictly) dominates the minimum risk equivariant predictive density estimator $\hat{q}_{\hbox{mre}}$.
\end{theorem}
{\bf Proof.}   Making use of (\ref{with-m}) and (\ref{J_1}) for the representation of $\hat{q}_{\pi_{U,A}}$, the  
difference in risks is given by
\begin{eqnarray*}
\Delta(\theta) &= & R_{KL}(\theta, \hat{q}_{\hbox{mre}}) - R_{KL}(\theta, \hat{q}_{\pi_{U,A}}) \\
\, & = &   \mathbb{E}^{X,Y_1} \log \left( \frac{\hat{q}_{\pi_{U,A}}(Y_1;X)}{\hat{q}_{\hbox{mre}}(Y_1;X)} \right) \\
\, & = &  \mathbb{E}^{X,Y_1} \log \left( \Phi(\alpha_0 + \alpha_1 \frac{Y_1 -X_1}{\tau}) - \Phi(\alpha_2 + \alpha_1 \frac{Y_1 -X_1}{\tau})\right) \\
\, & - &  \mathbb{E}^{X,Y_1} \log \left( \Phi(\frac{\alpha_0}{\sqrt{1+ \alpha_1^2}}) - \Phi(\frac{\alpha_2}{\sqrt{1+ \alpha_1^2}}) \right)\,,
\end{eqnarray*}
with the $\alpha_i$'s given in Section 2.2.  Now, observe that
\begin{equation}
\alpha_0 + \alpha_1 \frac{Y_1 -X_1}{\tau} = \frac{1}{\sigma_T} (m+ (X_1-X_2) + \beta (Y_1-X_1) \sim N(\delta_0= \frac{m+\theta_1-\theta_2}{\sigma_T}, 1)\,,
\end{equation}
and 
\begin{equation}
\frac{\alpha_0}{\sqrt{1+ \alpha_1^2}} = \frac{(m+ (X_1-X_2)}{\sqrt{\sigma_1^2 + \sigma_2^2}} \sim N(\delta_0'= \frac{m+\theta_1-\theta_2}{\sqrt{\sigma_1^2 + \sigma_2^2}}, 1)\,.
\end{equation}
Similarly,  we have $\alpha_2 + \alpha_1 \frac{Y_1 -X_1}{\tau} \sim N(\delta_2= \frac{-m+\theta_1-\theta_2}{\sigma_T}, 1)$ and $\frac{\alpha_2}{\sqrt{1+ \alpha_1^2}} \sim N(\delta_2'= \frac{-m+\theta_1-\theta_2}{\sqrt{\sigma_1^2 + \sigma_2^2}}, 1) $.  We thus can write  $$\Delta(\theta)\,=\, \mathbb{E} H(Z)\,,$$
$$ \hbox{ with }  H(Z)=  \log \left(\Phi(Z+\delta_0) - \Phi(Z+\delta_2) \right) \, - \, \log \left(\Phi(Z+\delta_0') - \Phi(Z+\delta_2') \right)\,,Z \sim N(0,1)\,.$$  
With $-m \leq \theta_1-\theta_2 \leq m$ and $\sigma_T^2 < \sigma_1^2 + \sigma_2^2$, we infer that $\delta_0 \geq \delta_0'$ with equality iff $\theta_1-\theta_2=-m$ and $\delta_2 \leq \delta_2'$ with equality iff $\theta_1-\theta_2=m$, so that $\mathbb{P}_{\theta}(H(Z)>0)=1$ and $\Delta(\theta)>0$ for all $\theta$ such that $|\theta_1-\theta_2|\leq m$. 
\qed

We now investigate situations where the variances in model (\ref{model}) are misspecified.  To this end, we consider $\sigma_1^2$, $\sigma_2^2$ and $\sigma_Y^2$ as the nominal variances used to construct the predictive density estimates  $\hat{q}_{\pi_{U,A}}$ and $\hat{q}_{\hbox{mre}}$, while the true variances, used to assess frequentist Kullback-Leibler risk, are, unbeknownst   to the investigator, given by 
$a_1^2 \sigma_1^2$, $a_2^2 \sigma_2^2$ and $a_Y^2 \sigma_Y^2$ respectively.  We exhibit, below in Theorem \ref{persistence}, many combinations of the nominal and true variances such that the Theorem \ref{Dominancecon}'s dominance result persists.  Such conditions for the dominance to persist includes the case of equal $a_1^2$, $a_2^2$ and $a_Y^2$ (i.e., the three ratios true variance over nominal variance are the same), among others.

We require the following intermediate result.

\begin{lemma} \label{U}
Let $U \sim N(\mu_U, \sigma_U^2)$ and $V \sim N(\mu_V, \sigma_V^2)$ with $\mu_U \geq \mu_V$  and $\sigma_U^2 \leq \sigma_V^2$.  Let $H$ be a differentiable function such that both $H$ and $-H'$ are increasing.  Then, we have $\mathbb{E} H(U) \geq \mathbb{E} H(V)$.
\end{lemma}
{\bf Proof.}  Suppose without loss of generality that $\mu_V=0$, and set  $s=\frac{\sigma_U}{\sigma_V}$.   Since $U$ and $\mu_U+sV$ share the same distribution and $\mu_U \geq 0$, we have:
\begin{eqnarray*}
\mathbb{E} H(U) & = & \mathbb{E} H(\mu_U+sV) \\
 \, & \geq & \mathbb{E} H(sV)\\
\,  & = &  \int_{\mathbb{R}_+}  \left( H(sv) + H(-sv)  \right)  \, \frac{1}{\sigma_V}  \phi(\frac{v}{\sigma_V})  \, dv \,. 
\end{eqnarray*}
 Differentiating with respect to $s$, we obtain
$$  \frac{d}{ds} \;  \mathbb{E} H(sV) \,=\, \int_{\mathbb{R}_+}  v \, \left( H'(sv) - H'(-sv)  \right)  \, \frac{1}{\sigma_V}  \phi(\frac{v}{\sigma_V})  \, dv \, \leq 0\,$$ 
since $H'$ is decreasing.  We thus conclude that
$$  \mathbb{E} H(U)  \geq \mathbb{E} H(sV) \geq \mathbb{E} H(V)\,,  $$
since $s\leq1$ and $H$ is increasing by assumption.  
\qed

\begin{theorem}
\label{persistence}
Consider model (\ref{model}) with $p=1$ and $A=[0,\infty)$.  Suppose that the variances are misspecified and that the true variances are given by  
$\mathbb{V}(X_1)= a_1^2 \sigma_1^2, \mathbb{V}(X_2)= a_2^2 \sigma_2^2, \mathbb{V}(Y_1)= a_Y^2 \sigma_Y^2$.
For estimating the density of $Y_1$ under Kullback-Leibler loss, the Bayes predictive density $\hat{q}_{\pi_{U,A}}$ dominates the minimum risk equivariant predictive density estimator $\hat{q}_{\hbox{mre}}$ whenever  $\sigma_U^2 \leq \sigma_V^2$ with
\begin{equation}
\label{persistencecondition}
\sigma_U^2 = \frac{a_2^2 \sigma_2^2 + (1-\beta)^2 a_1^2 \sigma_1^2 + \beta^2 a_Y^2 \sigma_Y^2}{\sigma_2^2+ \beta \sigma_Y^2}\,,\, \sigma_V^2= \frac{a_1^2 \sigma_1^2 + a_2^2 \sigma_2^2}{\sigma_1^2 + \sigma_2^2}\,, \beta= \frac{\sigma_1^2}{\sigma_1^2+\sigma_Y^2}. 
\end{equation}
In particular, dominance occurs for cases : {\bf (i)} $a_1^2=a_2^2=a_Y^2$, {\bf (ii)} $a_Y^2 \leq a_1^2=a_2^2 $, {\bf (iii)}  $\sigma_1^2=\sigma_2^2=\sigma_Y^2$ and $\frac{a_2^2+a_Y^2}{2} \leq a_1^2 $.  
\end{theorem}
\begin{remark}
Conditions {\bf (i)}, {\bf (ii)} and {\bf (iii)} are quite informative.  One common factor for the dominance to persist, especially seen by {\bf (iii)}, is for the variance of $X_1$ to be relatively large compared to the variances of $X_2$ and $Y_1$.
\end{remark}

{\bf Proof.}  Particular cases {\bf (i)}, {\bf (ii)}, {\bf (iii)} follow easily from (\ref{persistencecondition}).   To establish condition (\ref{persistencecondition}),  we prove, as in Theorem \ref{Dominancecon}, that $\Delta(\theta)$ given in (\ref{diffrisks}) is greater or equal to zero.  We apply Lemma \ref{U}, with $H \equiv \log \Phi$ increasing and concave  as required, showing that $\mathbb{E} \, log \Phi(U) \geq \mathbb{E} \, log \Phi(V)$ with $U= \alpha_0 + \alpha_1 \frac{Y_1 -X_1}{\tau} \sim N(\mu_U, \sigma_U^2)$ and $V= \frac{\alpha_0}{\sqrt{1+\alpha_1^2}} \sim N(\mu_V, \sigma_V^2)$.  Since $\mu_U = \frac{\theta_1-\theta_2}{\sigma_T} > \frac{\theta_1-\theta_2}{\sqrt{\sigma_1^2+\sigma_2^2}} = \mu_V $, the inequality $\sigma_U^2 \leq \sigma_V^2$ will suffice to have dominance.  Finally, the proof is complete by checking that $\sigma_U^2$ and $\sigma_V^2$ are as given in (\ref{persistencecondition}), when the true variances are given by
$\mathbb{V}(X_1)= a_1^2 \sigma_1^2, \mathbb{V}(X_2)= a_2^2 \sigma_2^2, \mathbb{V}(Y_1)= a_Y^2 \sigma_Y^2$.   \qed

\begin{remark}
In opposition to the above robustness analysis, the dominance property of $\hat{q}_{\pi_{U,A}}$ versus $\hat{q}_{\hbox{mre}}$ for the restriction $\theta_1-\theta_2 \geq 0$ does not persists for parameter space values such that $\theta_1-\theta_2<0$, i.e., the additional information difference is misspecified.  In fact, it is easy to see following the proof of Theorem \ref{Dominancecon} that $R_{KL}(\theta, \hat{q}_{\hbox{mre}}) - R_{KL}(\theta, \hat{q}_{\pi_{U,A}})<0$ for $\theta$'s such that $\theta_1-\theta_2<0$.
A potential protection is to use the predictive density estimator $\hat{q}_{\pi_{U,A'}}$ with $A'=[\epsilon,\infty)$, $\epsilon<0$, and with dominance occurring for all $\theta$ such that $\theta_1-\theta_2 \geq \epsilon$ (Remark \ref{transformation} and Theorem \ref{Dominancecon}).  
\end{remark}

\section{Examples, illustrations and further comments}

We present and comment numerical evaluations of Kullback-Leibler risks in the univariate case for both $\theta_1 \geq \theta_2$ (Figures 1, 2) and $|\theta_1 - \theta_2| \leq m, m=1,2.$ (Figures 3, 4).  Each of the figures consists of plots of risk ratios, as functions of $\Delta=\theta_1 - \theta_2$ with the benchmark $\hat{q}_{mre}$ as the reference point.  The variances are set equal to $1$, except for Figure 2 which highlights the effect of varying $\sigma^2_2$.

Figure 1 illustrates the effectiveness of variance expansion (Corollary \ref{cormledominance}), as well as the dominance finding of Theorem \ref{Dominancecon}.  More 
precisely, the Figure relates to Example \ref{example-univariatecase} where $\hat{q}_{mle}$ is improved by the variance expansion version $\hat{q}_{mle, 2}$, which belongs both to the subclass of dominating densities $\hat{q}_{mle, c}$ as well as to the complete subclass of such predictive densities.  The gains are impressive ranging from a minimum of about $8\%$ at $\Delta=0$ to a supremum value of about $44\%$ for $\Delta \to \infty$.  Moreover, the predictive density $\hat{q}_{mle, 2}$ also dominates $\hat{q}_{mre}$ by duality, but the gains are more modest.  Interestingly, the penalty of failing to expand is more severe than the penalty for using an inefficient {\em plug-in} estimator of the mean.  In accordance with Theorem \ref{Dominancecon}, the Bayes predictive density $\hat{q}_{\pi_{U,A}}$ improves uniformly on $\hat{q}_{mre}$ except at $\Delta=0$ where the risks are equal.  As well, $\hat{q}_{\pi_{U,A}}$ compares well to $\hat{q}_{mle, 2}$, except for small $\Delta$,  with
$R(\theta, \hat{q}_{mle, 2}) \leq R(\theta, \hat{q}_{\pi_{U,A}})$  if and only if $\Delta \leq \Delta_0$ with $\Delta_0 \approx 0.76$.

Figure 2 compares the efficiency of the predictive densities $\hat {q}_{\pi_{U,A}}$ and $\hat{q}_{mre}$ for varying $\sigma^2_2$.  Smaller values of $\sigma^2_2$ represent more precise estimation of $\theta_2$ and translates to a tendency for the gains offered by $\hat {q}_{\pi_{U,A}}$ to be greater for smaller $\sigma^2_2$; but the situation is slightly reversed for larger $\Delta$. 

\begin{figure}[H] 
    \centering 
    \includegraphics[width=0.7\textwidth]{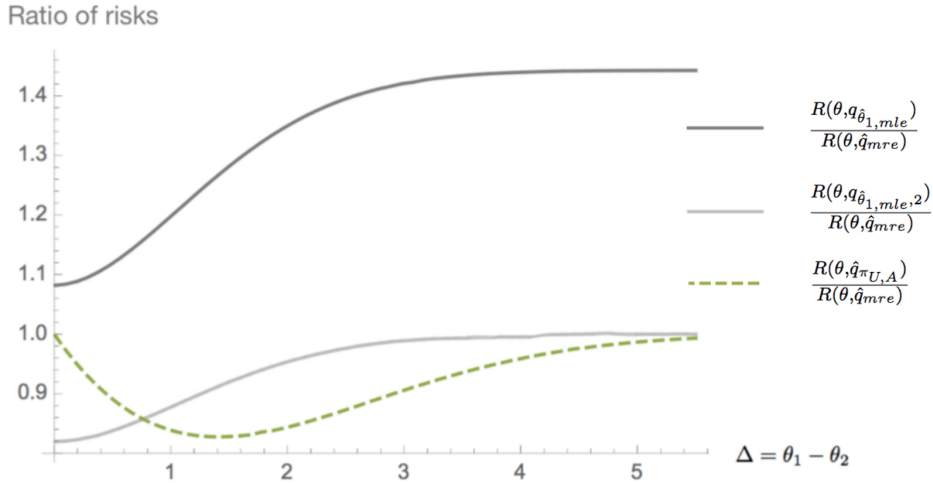}
    \caption{Kullback-Leibler risk ratios for $p=1$, $A=[0,\infty)$, and $\sigma^2_1=\sigma^2_2=\sigma^2_Y=1$}
    \label{Figure1}  
\end{figure}

\begin{figure}[H] 
   \centering  
   \includegraphics[width=0.7\textwidth]{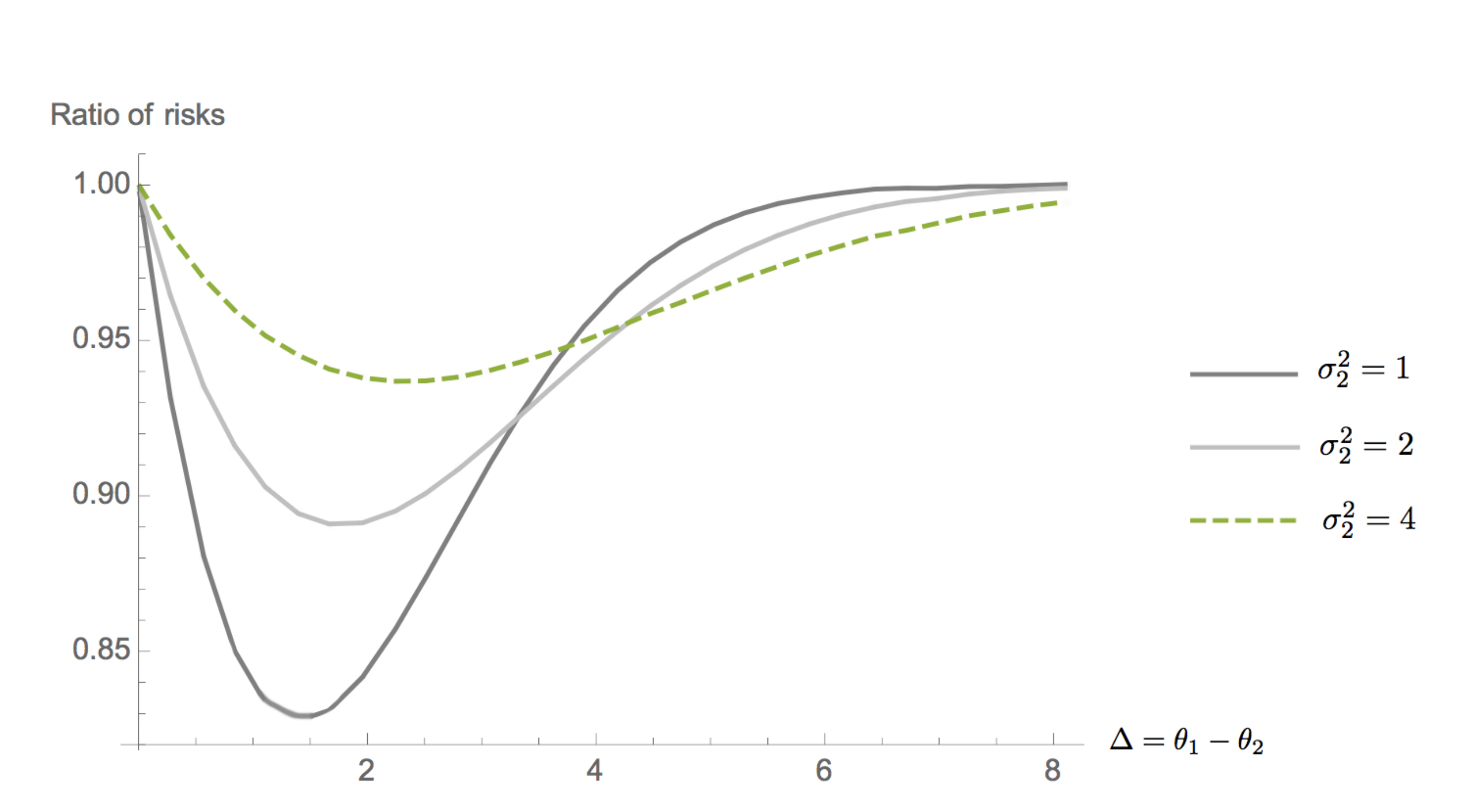}
    \caption{Kullback-Leibler risk ratios for $p=1$, $A=[0,\infty)$, $\sigma^2_1=\sigma^2_Y=1$ and $\sigma^2_2=1,2,4$ } \label{Figure2}      
\end{figure}

Figures 3 and 4 compare the same estimators as in Figure 1, but they are adapted to the restriction to compact interval.  Several of the features of Figure 1 are reproduced with the noticeable inefficiency of $\hat{q}_{mle}$ compared to both $\hat{q}_{mle, 2}$ and $\hat{q}_{\pi_{U,A}}$.  For the larger parameter space (i.e. $m=2$), even $\hat{q}_{mre}$ outperforms $\hat{q}_{mle}$ as in Figure 1, but the situation is reversed for $m=1$ where the efficiency of  better point maximum likelihood estimates plays a more important role.  The Bayes performs well, dominating $\hat{q}_{mre}$ in accordance with Theorem \ref{dominancedb}, especially for small of moderate $\Delta$, and even improving on $\hat{q}_{mle, 2}$ for $m=1$.  Finally, we have extended the plots outside the parameter space which is useful for assessing performance for slightly incorrect specifications of the additional information.

\begin{figure}[H] 
   \centering  
    \includegraphics[width=0.7\textwidth]{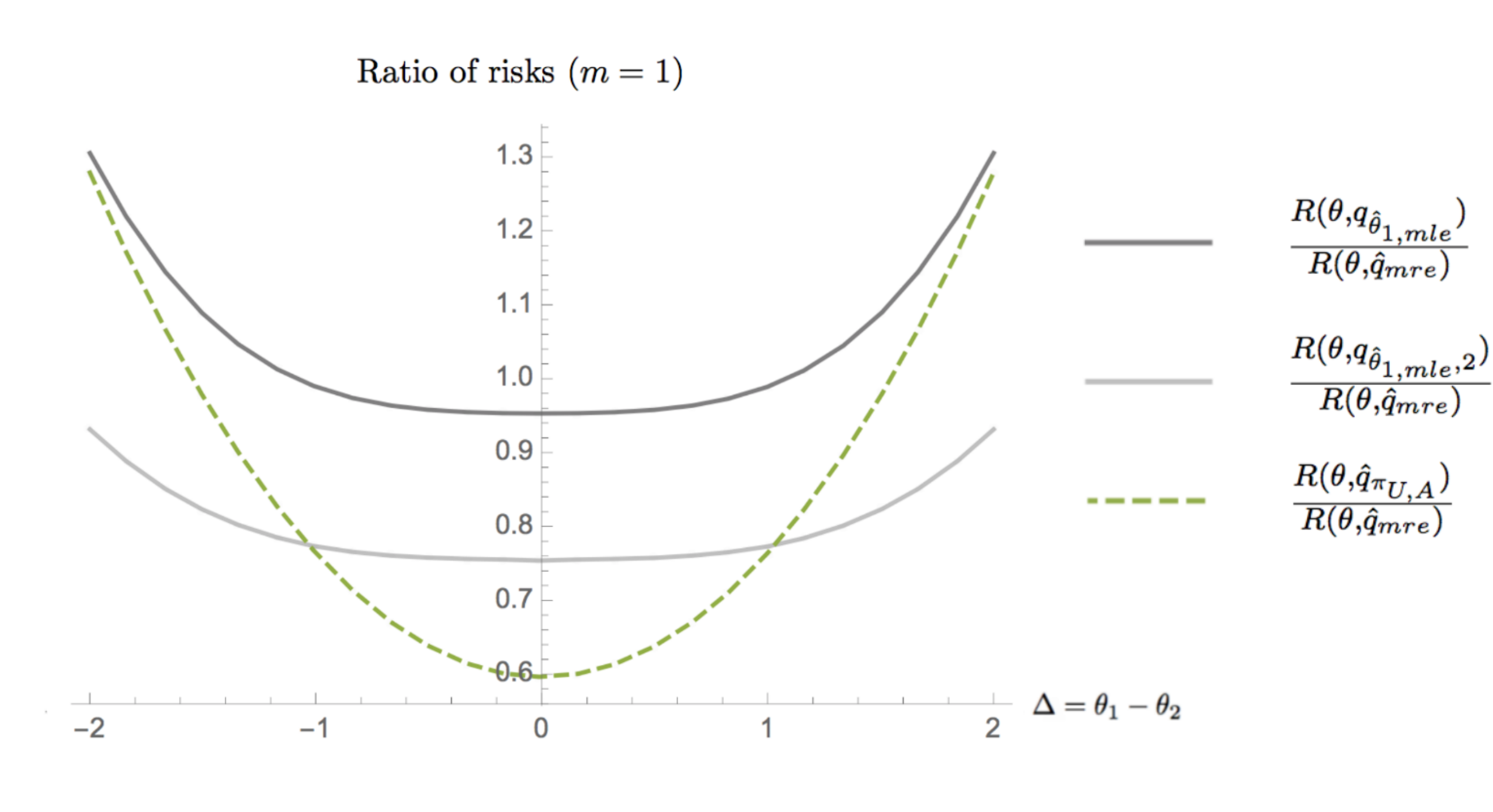}
    \caption{Kullback-Leibler risk ratios for $p=1$, $A=[-1,1]$, and $\sigma^2_1=\sigma^2_2=\sigma^2_Y=1$} \label{Figure3}      
\end{figure}

\begin{figure}[H] 
   \centering  
   \includegraphics[width=0.8\textwidth]{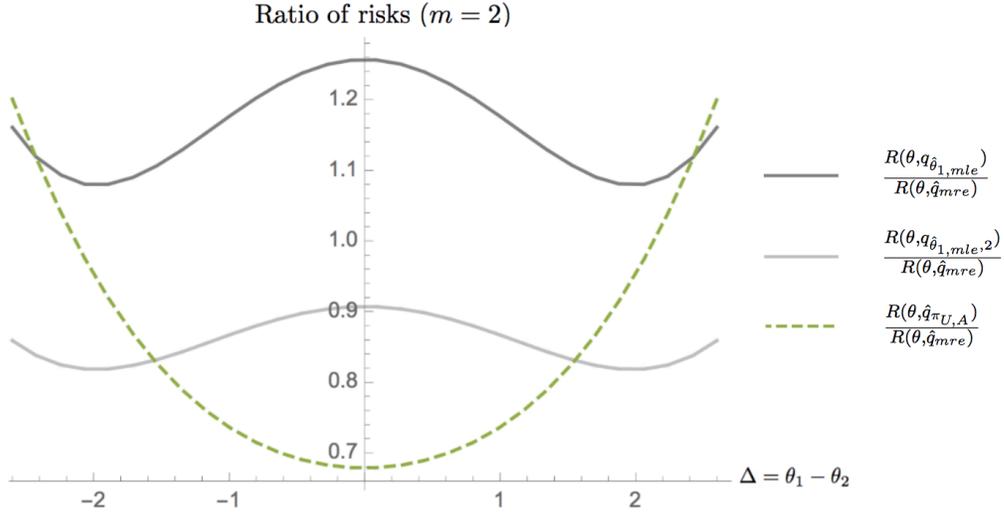}
    \caption{Kullback-Leibler risk ratios for $p=1$, $A=[-2,2]$, and $\sigma^2_1=\sigma^2_2=\sigma^2_Y=1$ } \label{Figure4}      
\end{figure}

\section{Concluding remarks}

For multivariate normal observables $X_1 \sim N_p(\theta_1, \sigma^2_1 I_p)$, $X_2 \sim N_p(\theta_2, \sigma^2_2 I_p)$, we have provided findings concerning the efficiency of predictive density estimators $Y_1 \sim N_p(\theta_1, \sigma^2_1 I_p)$ with the added parametric information $\theta_1 - \theta_2 \in A$.  Several findings provide improvements on benchmark predictive densities, such those obtained as {\em plug-in's}, as maximum likelihood, or as minimum risk equivariant.  The results range over a class of $\alpha-$divergence losses, different settings for $A$, and include Bayesian improvements for reverse Kullback-Leibler and Kullback-Leibler losses.  The various techniques used lead to novel connections between different problems, which is also of interest as, for instance,  described following both Proposition \ref{propositionKL-RKL} and Proposition \ref{propositionalpha}-Remark \ref{remark |alpha|=1}. 

Although the Bayesian dominance results for Kullback-Leibler loss for $p=1$ extend to the rectangular case with $\theta_{1,i}-\theta_{2,i} \in A_i$ for $i=1,\ldots, p$ and the $A_i's$ either lower bounded, upper bounded, or bounded to intervals $[-m_i,m_i]$ 
(since the Kullback-Leibler divergence for the joint density of $Y$ factors and becomes the sum of the marginal Kullback-Leibler divergences, and that the posterior distributions of the $\theta_{1,i}$'s are independent), a general Bayesian dominance result of $\hat{q}_{\pi_{U,A}}$ over $\hat{q}_{mre}$, is lacking and would be of interest.  As well,  comparisons of predictive densities for the case of homogeneous, but unknown variance (i.e., $\sigma^2_1=\sigma^2_2=\sigma^2_Y$), are equally of interest.  Finally, the analyses carried out here should be useful as benchmarks in situations where the constraint set $A$ has an anticipated form, but yet is unknown.  In such situations, a reasonable approach would be to consider priors that incorporate uncertainty on $A$, such as setting $A=\{\theta: \|\theta_1 - \theta_2 \| \leq m\}$, $A=[m, \infty)$, with prior uncertainty specified for $m$.

\section*{Appendix}
\subsection*{Predictive density estimation under reverse Kullback-Leibler loss}

\noindent  The objective of this part is two-fold.  First, we present a quite general result which stipulates that Bayes predictive density estimators are always \emph{plug-in} densities in an exponential family set-up with, or without, additional information.  Such a result was obtained by Yanagimoto and Ohnishi (2009).  We provide an extension for problems with additional information and we seek to give more prominence to Yanagimoto and Ohnishi's wonderful result.  
Secondly, applications of Theorems \ref{Bayes=plugin} and \ref{Bayes=plugin-D} yield part {\bf (b)} of Lemma \ref{bayesgeneral} and 
the reverse Kullback-Leibler part of Lemma \ref{klrkldual}.

\noindent Consider the exponential family model densities, with respect to $\sigma-$finite measures $\mu_1$ and $\mu_2$, under canonical form
\begin{eqnarray}
\label{modelexp}
\nonumber X|\eta &\sim& p_{\eta}(x)=h_1(x) \exp\{{\eta_1}^T s_1(x) + {\eta_2}^T s_2(x) - c_1(\eta) \},  \vspace{2cm}\\
Y_1|\eta_1 &\sim& q_{\eta_1}(y)=h_2(y_1) \exp\{ {\eta_1}^T t_1(y_1) - c_2(\eta_1) \}\,, 
\end{eqnarray}
where $X=(X_1,X_2)^T$, $\eta=(\eta_1,\eta_2)^T$, and $\eta_1,\, \eta_2,\, s_1(x),\, s_2(x), \,t_1(y_1)$ are vectors of dimension $p$.  
In this set-up, we assume that $X$ and $Y_1$ are independently distributed given $\eta$, $\eta_1$ is a common parameter, and we seek a predictive density for $Y_1$ based on $X$ and with the additional information $\eta_1-\eta_2 \in A$. 
We thus consider predictive densities $\hat{q}(\cdot;X)$ for $Y_1$
and their performance as evaluated by reverse Kullback-Leibler loss 
\begin{equation}
\label{lk}
L(\eta_1,\hat{q}) = \int \hat{q}(y_1) \log\left(\frac{\hat{q}(y_1)}{q_{\eta_1}(y_1)}  \right) \, d\mu_2(y_1)\,,
\end{equation} 
and corresponding risk
$$ R(\eta, \hat{q}) = \int \int p_{\eta}(x) \hat{q}(y_1;x) \log\left(\frac{\hat{q}(y_1;x)}{q_{\eta_1}(y_1)}  \right) \, d\mu_2(y_1)\,d\mu_1(x)\,. $$
A \emph{ plug--in} estimator for the density $q_{\eta_1}$ is simply of the form $q_{\hat{\eta}_1(X)}$.  For Kullback-Leibler loss, obtained by switching $q_{\eta_1}$ and $\hat{q}$  in (\ref{lk}), \emph{plug-in} density estimators are not compatible with Bayesianity and can be quite inefficient in terms of Kullback-Leibler risk, as seen above in Lemma \ref{c} for normal models.  However, for reverse Kullback-Leibler loss, the situation is the opposite, and universally so for the exponential family set-up above as shown in Theorem \ref{bayesrepresentationrkl}.  Furthermore, the \emph{ plug--in} estimator is the posterior expectation of $\eta_1$.   This holds regardless of the prior on $\eta$ (including cases where $\eta_1-\eta_2 \in A$) and the particular forms of $p_{\eta}$ and $q_{\eta_1}$.  This was observed and exploited for normal models by Maruyama and Strawderman (2012). 

The second observation made below concerns the frequentist risk of \emph{ plug--in} densities.  Indeed, reverse Kullback-Leibler loss (among others) for a \emph{plug-in} estimate becomes simply a measure of distance between the densities $q_{\eta_1}$ and $q_{\hat{\eta_1}}$, otherwise known as intrinsic loss (e.g. Robert, 1996).  For exponential families, as noted by Brown (1986, Proposition 6.3),  such a distance has a simple and appealing form.  Here, it leads to a representation, for both \emph{ plug--in} and thus Bayes predictive density estimators,  of the reverse Kullback-Leibler risk in terms of the point estimate risk performance of the same \emph{plug--in} estimator with respect to a dual loss.    

The following representation of a Bayes predictive density estimator under reverse Kullback-Leibler is well known (e.g., Corcuera and Giummol\`e, 1999), but we provide a short presentation for completeness.

\begin{lemma}
\label{bayesrepresentationrkl}
For estimating $q_{\eta_1}$ under reverse Kullback-Leibler loss and based on $X$ as in (\ref{modelexp}), the Bayes predictive density estimator is $ \hat{q}_{\pi}(y_1;x) \propto \exp{\{E(\log q_{\eta_1}(y_1)|x) \}}\,$. 
\end{lemma}
{\bf Proof.} For an estimator $\hat{q}$ and denoting $G_x$ as the posterior c.d.f. of $\eta$, the expected posterior loss may be expressed as:
\begin{eqnarray}
\nonumber E \left(L(\eta_1, \hat{q})|x\right) & = & \int \{\int  \hat{q}(y_1) \left(\log \hat{q}(y_1) - \log q_{\eta_1}(y_1)\right) \, d\mu_2(y_1) \,\}  \, dG_x(\eta) \,  \\
\nonumber \, &=& \int \hat{q}(y_1) \, \{\, \log \hat{q}(y_1) \, - \,E(\log q_{\eta_1}(y_1)|x)\, \}  \, d\mu_2(y_1)  \\
\label{last} \, &=&  \log c + \int \hat{q}(y_1) \, \{- \log (\frac{\hat{q}_{\pi}(y_1;x)}{\hat{q}(y_1)})  \}  \, d\mu_2(y_1) \,,
\end{eqnarray}
where $ \hat{q}_{\pi}(y_1;x) =c \,\exp{\{E(\log q_{\eta_1}(y_1)|x) \}}\,$.
Using Jensen's inequality applied to $-\log$, we obtain indeed from (\ref{last}), for all estimators
$\hat{q}$, 
$$ E \left( L(\eta_1, \hat{q}\,)|x\right) \, \geq \, \log \,\,c \, - \, \log \int \hat{q}_{\pi}(y_1;x) \, d\mu_2(y_1) = \,\log \, c \, = \, E \left( L(\eta_1, \hat{q}_{\pi})|x\right)\,.  \qed $$

The following representation applies with or without the additional information provided by the constraint $\eta_1 - \eta_2 \in A$, with the additional information case representing an extension of Yanagimoto and Ohnishi's result.  

\begin{theorem}
\label{Bayes=plugin}
For model (\ref{modelexp}), reverse Kullback-Leibler loss,  a prior measure $\pi$ for $\eta$ such that the posterior distribution and expectation exists, the Bayes predictive density estimate $\hat{q}_{\pi}(\cdot;x)$ is the \emph{ plug--in} density estimate $q_{\hat{\eta}_1}(\cdot;x)$, with $\hat{\eta}_1(x)=E_{\pi}(\eta_1|x)$ the posterior expectation of $\eta_1$.
\end{theorem}
{\bf Proof.}  Using Lemma \ref{bayesrepresentationrkl}, we obtain
\begin{eqnarray*}
\hat{q}_{\pi}(y_1;x) &\propto& \exp{\{E(\log q_{\eta_1}(y_1)|x) \}} \\
\, & \propto & h_2(y_1) \exp\{E(\eta_1^T t_1(y_1) - c_2(\eta_1)|x)\} \\
\, & \propto & h_2(y_1)  \exp\{(\eta_1^T t_1(y_1) - c_2(E(\eta_1)|x)) \},
\end{eqnarray*}
which matches indeed the \emph{plug-in} density $q_{\hat{\eta}_1}(\cdot;x)$ with $\hat{\eta}_1(x)=E_{\pi}(\eta_1|x)$.  \qed

\begin{theorem}
\label{Bayes=plugin-D}
For model (\ref{modelexp}), the reverse Kullback-Leibler frequentist risk of the \emph{ plug--in} density $q_{\hat{\eta}_1}(\cdot;X)$ is equivalent to the frequentist risk for estimating $\eta_1$ based on $X$ under the dual point estimation loss
$$
L_{dual}(\eta_1, \hat{\eta}_1)=(\hat{\eta}_1-\eta_1)^T \,\mathbb{E}_{\hat{\eta}_1} (t(Y))+(c_2(\eta_1)-c_2(\hat{\eta}_1)).
$$ 
\end{theorem}

{\bf Proof.}
For the plug--in density estimator, we have 
\begin{eqnarray*}
L_{dual}(\eta_1, \hat{\eta}_1)&=&\int q_{\hat{\eta}_1}(y_1) \log \frac{q_{\hat{\eta}_1}(y_1)}{q_{\eta_1}(y_1)} \,d \mu_2(y_1)\\
&=& \int q_{\hat{\eta}_1}(y_1) \{(\hat{\eta}_1-\eta_1)^T\, t(y_1)+(c_2(\eta_1)-c_2(\hat{\eta}_1))\} \, \,d \mu_2(y_1) \\
&=& (\hat{\eta}_1-\eta_1)^T\, \mathbb{E}_{\hat{\eta}_1} t(Y_1)+(c_2(\eta_1)-c_2(\hat{\eta}_1))\,,
\end{eqnarray*}
which leads to the result.  \qed

\begin{example}		
\label{exampleappendix}
For the multivariate normal model (\ref{model}), the last two theorems apply as examples of model (\ref{modelexp}) with $\eta_1=\theta_1$, $\eta_2=\theta_2$, $c_2(\eta_1)=\frac{\|\eta_1\|^2}{2 \sigma^2_Y}$, $t(y_1)=\frac{y_1}{\sigma^2_Y}$.  Theorem \ref{Bayes=plugin} yields the Bayes predictive density given in (\ref{bayesrkl}), while 
Theorem \ref{Bayes=plugin-D} yields the dual loss
 $L_{dual}(\eta_1, \hat{\eta}_1)=  (\hat{\eta}_1 - \eta_1)^T \, \mathbb{E}_{\hat{\eta}_1} (\frac{Y_1}{\sigma^2_Y}) + \frac{\|\eta_1\|^2}{2 \sigma^2_Y} - \frac{\|\hat{\eta}_1\|^2}{2 \sigma^2_Y}  =  \frac{\norm{\hat{\eta}_1-\eta_1}}{2\sigma^2}$, 
as stated in Lemma \ref{klrkldual}. 
\end{example}

\section*{Acknowledgments}
Author Marchand gratefully acknowledges the research support from the Natural Sciences and Engineering Research Council of Canada.   We are grateful to Bill Strawderman for useful discussions, namely on the developments for reverse Kullback-Leibler loss in the Appendix.

\bibliographystyle{plain}
\bibliography{dissertation}

\end{document}